\documentclass[
               DIV=15,
               ]{scrartcl}  

\usepackage{color} 
\usepackage{psfrag} 
\usepackage{subcaption}
\usepackage{hyperref}
\usepackage{amsmath}
\usepackage{amssymb}
\usepackage{amsthm}
\usepackage{graphics}
\usepackage{todonotes}
\usepackage{amsfonts}  
\usepackage{adjustbox}    
\usepackage{mathtools}
\usepackage{xcolor}        
\usepackage{blkarray}        
\usepackage{footnote}
\makesavenoteenv{tabular}
\makesavenoteenv{table}

\usepackage[square,numbers,comma]{natbib}

\hypersetup{colorlinks,citecolor=blue}

\usepackage{standalone}
 
\usepackage{pgfplots}
\usepackage{pgfplotstable}
\usetikzlibrary{pgfplots.groupplots}
\usetikzlibrary{fit}
\usepackage{filecontents}
\usepackage{tikz-3dplot}

\usepgfplotslibrary{patchplots}

\usepackage{tikz}
\usepackage{tikzscale}

\usetikzlibrary{arrows,3d}

\usetikzlibrary{shapes.geometric}
\usetikzlibrary{positioning}
\usetikzlibrary{calc}
\usetikzlibrary{decorations.pathreplacing} 
\usetikzlibrary{decorations.markings}

\usetikzlibrary{spy}
\usetikzlibrary{snakes,arrows,shapes}    
\usetikzlibrary{spy}
\usetikzlibrary{backgrounds,tikzmark}  
\usetikzlibrary{decorations}
\usetikzlibrary{positioning}
\usetikzlibrary{patterns}
\usetikzlibrary{calc} 

\usepackage{booktabs}   
\usepackage{makecell} 
\usepackage{colortbl}   

\usepackage{xspace}
\usepackage{color}     
\usepackage{setspace}   
 
\usepackage{soul}
\usepackage{ulem}
\usepackage{currfile}
\usepackage{import}

\usepackage{authoraftertitle} 
\usepackage{authblk} 

\usepackage[noabbrev]{cleveref}

\usepackage{enumitem}
\usepackage{multicol} 
\usepackage{algorithm}
\usepackage{algpseudocode}
\usepackage{algorithmicx}
\usepackage[section]{placeins}

\usepackage{environ}

    \usepackage{framed}
    
    \setlist{noitemsep,topsep=2pt,parsep=0pt,partopsep=0pt}
  
    \usepackage{lmodern}

    \usetikzlibrary{external}
    \tikzset{external/force remake}
    
   \pgfdeclarelayer{background}
   \pgfdeclarelayer{foreground}
   \pgfsetlayers{background,main,foreground}

\pgfplotsset{compat=1.8}

\newcommand{\findmax}[3]{
    \pgfplotstablesort[sort key={#2},sort cmp={float >}]{\sorted}{#1}%
    \pgfplotstablegetelem{0}{#2}\of{\sorted}%
    \let #3=\pgfplotsretval%
}

\pgfplotscreateplotcyclelist{myCycleListBlackAndWhite}
{%
  {black, mark=o},
  {black, mark=square},
  {black, mark=triangle},
  {black, mark=diamond}, 
  {black, mark=pentagon}, 
  {black, mark=*},
  {black, mark=square*},
  {black, mark=triangle*},
  {black, mark=diamond*}, 
  {black, mark=pentagon*}, 
}

\definecolor{darkgreen}{rgb}{0,0.4,0} 
\definecolor{darkbrown}{rgb}{0.5, 0.396, 0.09}

\definecolor{c1}{rgb}{0.0, 0.4196078431372549, 0.6431372549019608}
\definecolor{c2}{rgb}{1.0, 0.5019607843137255, 0.054901960784313725}
\definecolor{c3}{rgb}{0.6705882352941176, 0.6705882352941176,
0.6705882352941176} \definecolor{c}{rgb}{0.34901960784313724, 0.34901960784313724, 0.34901960784313724}
\definecolor{c4}{rgb}{0.37254901960784315, 0.6196078431372549,
0.8196078431372549} \definecolor{c}{rgb}{0.7843137254901961, 0.3215686274509804, 0.0}
\definecolor{c5}{rgb}{0.5372549019607843, 0.5372549019607843,
0.5372549019607843} \definecolor{c}{rgb}{0.6352941176470588, 0.7843137254901961, 0.9254901960784314}
\definecolor{c6}{rgb}{1.0, 0.7372549019607844, 0.4745098039215686}
\definecolor{c7}{rgb}{0.8117647058823529, 0.8117647058823529,
0.8117647058823529}

\pgfplotscreateplotcyclelist{myCycleListColor}
{%
  {blue, mark=o},
  {red, mark=square},
  {darkbrown, mark=triangle},
  {darkgreen, mark=diamond},
  {black, mark=pentagon},
  {blue, mark=*},
  {red, mark=square*},
  {darkbrown, mark=triangle*},
  {darkgreen, mark=diamond*},
  {black, mark=pentagon*},
}

\pgfplotsset{every axis/.append style= 
              {
                font=\scriptsize,
                mark size=3,
                legend style={font=\tiny, mark size=3, draw=none, fill=none},
                legend cell align=left,
                cycle list name=myCycleListColor,
                scaled y ticks = false,
				scaled x ticks = false,
				trim axis left,
				trim axis right,
				sharp plot,
				tick label style ={font=\tiny},
				label style ={font=\scriptsize},
				very thin,
				ymajorgrids=true,
				grid style=dotted,
				legend pos= north east,
				height=\figHeight, width=\figWidth,
              }
            }

\pgfplotsset{every tick label/.append style={font=\tiny}}
\pgfplotsset{every label/.append style={font=\scriptsize}}
\pgfplotsset{every axis legend/.append style={font=\tiny}}
\pgfplotsset{every axis plot/.append style={thick}}
            
\pgfplotstableset
{
  every odd row/.style={before row={\rowcolor[gray]{0.9}}},
  every head row/.style=
  {
    before row=
    {%
      \toprule
    },
    after row=\midrule
  },
  every last row/.style=
  {
    after row=\bottomrule
  }
}

\newif\ifdrawboundingbox

\pgfkeys
{ 
  /pgf/number format/.cd,
  fixed,
  set thousands separator={\,},
  1000 sep in fractionals,
}

\newcolumntype{C}[1]{>{\centering\arraybackslash}m{#1}}
\newcolumntype{R}[1]{>{\raggedright\arraybackslash}m{#1}}
\newcolumntype{L}[1]{>{\raggedleft\arraybackslash}m{#1}}

\newcommand\restr[2]{{
		\left.\kern-\nulldelimiterspace 
		#1 
		\right|_{#2} 
}}

\newcommand\compactDots{\ifmmode\ldots\else\makebox[10cm][c]{.\hfil.\hfil.}\fi}

\newcommand{\supp}[1]{\textnormal{supp}(#1)}

\def\registered{{\ooalign{\hfil\raise .00ex\hbox{\tiny R}\hfil\cr\mathhexbox20D}}}

\makeatletter
\def\namedlabel#1#2{\begingroup
	\def\@currentlabel{#2}%
	\label{#1}\endgroup
}
\makeatother

\makeatletter
\renewcommand{\todo}[2][]{\tikzexternaldisable\@todo[#1]{#2}\tikzexternalenable}
\makeatother

\usepackage{etoolbox}
\AtBeginEnvironment{mysmallmatrix}{%
\setstretch{0.9}%
\setlength{\arraycolsep}{2pt}%
}%

{%
	\begin{matrix}%
}%
	{\end{matrix}}%

{%
	\begin{matrix}%
	}%
	{\end{matrix}}%
\AtBeginEnvironment{myverysmallmatrix}{%
	\footnotesize%
	\setstretch{0.55}%
	\setlength{\arraycolsep}{1.5pt}%
}%

\newcommand{\overbar}[1]{\mkern 0.1mu\overline{\mkern-0.1mu#1\mkern-0.1mu}\mkern 0.1mu}
\makeatletter
\def\bbordermatrix#1{\begingroup \m@th
	\@tempdima 4.75\p@
	\setbox\z@\vbox{%
		\def\cr{\crcr\noalign{\kern2\p@\global\let\cr\endline}}%
		\ialign{$##$\hfil\kern2\p@\kern\@tempdima&\thinspace\hfil$##$\hfil
			&&\quad\hfil$##$\hfil\crcr
			\omit\strut\hfil\crcr\noalign{\kern-\baselineskip}%
			#1\crcr\omit\strut\cr}}%
	\setbox\tw@\vbox{\unvcopy\z@\global\setbox\@ne\lastbox}%
	\setbox\tw@\hbox{\unhbox\@ne\unskip\global\setbox\@ne\lastbox}%
	\setbox\tw@\hbox{$\kern\wd\@ne\kern-\@tempdima\left[\kern-\wd\@ne
		\global\setbox\@ne\vbox{\box\@ne\kern2\p@}%
		\vcenter{\kern-\ht\@ne\unvbox\z@\kern-\baselineskip}\,\right]$}%
	\null\;\vbox{\kern\ht\@ne\box\tw@}\endgroup}
\makeatother

\title{A hierarchical approach to the \textit{a posteriori} error estimation of isogeometric Kirchhoff plates and Kirchhoff-Love shells}

\author[1]{Pablo Antolin}
\author[1,2]{Annalisa Buffa}
\author[1]{Luca Coradello\thanks{luca.coradello@epfl.ch, Corresponding Author}}

\affil[1]{Institute of Mathematics,
 		  \'Ecole Polytechnique F\'ed\'erale de Lausanne, Lausanne, Switzerland.}

\affil[2]{Istituto di Matematica Applicata e Tecnologie Informatiche `E. Magenes' (CNR), Pavia, Italy.}

\newcommand{\publicationDate}{\today}

\usepackage{fancyhdr}
\date{}
\fancypagestyle{plain}{%
  \fancyhf{}%
  \fancyfoot[L]{
  \vspace{-1.25cm} 
  \footnotesize{Preprint submitted to \textit{}  \hfill
  \publicationDate
  \\
  }} }
 
\begin{document}  
\normalem
\maketitle  
  
\vspace{-1.5cm} 
\hrule 
\section*{Abstract}
This work focuses on the development of a posteriori error estimates for fourth-order, elliptic, partial differential equations. In particular, we propose a novel algorithm to steer an adaptive simulation in the context of Kirchhoff plates and Kirchhoff-Love shells by exploiting the local refinement capabilities of hierarchical B-splines. The method is based on the solution of an auxiliary residual-like variational problem, formulated by means of a space of localized spline functions. This space is characterized by $C^1$ continuous B-splines with compact support on each active element of the hierarchical mesh.   
We demonstrate the applicability of the proposed estimator to Kirchhoff plates and Kirchhoff-Love shells by studying several benchmark problems which exhibit both smooth and singular solutions. In all cases, we obtain optimal asymptotic rates of convergence for the error measured in the energy norm and an excellent approximation of the true error.
\vspace{1cm}

\noindent \textit{Keywords}: isogeometric analysis, a posteriori error estimator, adaptivity, hierarchical B-splines, shells, plates.
\vspace{0.2cm} 
\hrule

\def\Estconst{3}
 


\section{Introduction}\label{sec:introduction}

IsoGeometric Analysis (IGA) has been a thriving area of research since the first pioneering work \citep{Hughes2005} was published in 2005. Employing smooth B-splines, NURBS or variances thereof as basis functions for the solution field has shown excellent properties in many mathematical and engineering applications. For instance, in many cases higher global continuity yields a better accuracy per degree-of-freedom (dof) and it allows to discretize higher-order variational problems in their primal form, e.g. Kirchhoff plates and Kirchhoff-Love shells \citep{Kiendl2009} or Cahn-Hilliard problems \citep{Gomez2008}. For a detailed review of the method and its applications, the reader is referred to \citep{Hughes2005,Cottrell2009,Hughes2017special}.

However, in order to make IGA competitive for industrial applications, a crucial feature is the capability of building locally refined bases suitable for analysis in an efficient and robust manner. Indeed, we remark that the tensor-product nature of B-splines is broken by local refinement. To tackle these issues, several alternatives have been proposed such as hierarchical B-splines (HB) \citep{Forsey1988, Greiner1997, Kraft1997} and their recent variant denoted by truncated hierarchical B-splines (THB) \citep{Giannelli2012, Giannelli2016}, T-splines \citep{Bazilevs2010,Scott2012,Beirao2013} and locally-refinable splines (LR-splines) \citep{Dokken2013,Bressan2013}. Also the inverse operation of refinement (known as coarsening), which plays an important role for the computational efficiency of dynamic problems, proves to be a non-trivial task and another booming area of research, see for instance \citep{Lorenzo2017,Garau2018,Hennig2018,Carraturo2019}. In the framework of this work, we employ hierarchical splines as they provide the most straightforward construction and an easy implementation in the context of isogeometric analysis.

In order to automate the process of local refinement, \textit{a posteriori} error estimators have been successfully used to steer adaptive simulations for standard finite elements, where the literature is rich and several families of estimators have been thoroughly studied (see for instance \citep{Babuska1981,Ainsworth1997,verfurth2013} and references therein). We remark that a good estimator should properly capture and resolve local features of the solution that cause the error to be much larger than in the rest of the computational domain.   
In the realm of isogeometric analysis, the field of error estimation is still in its early stages of development. Adaptivity using T-splines has been presented in the pioneering work \citep{Juettler2010}. More results on recovery-based error estimators have been published in \citep{Kumar2015,Kumar2017} where refinement is achieved by using LR-splines. Additionally, functional-type error estimates have been studied in \citep{Kleiss2015} but not within a real adaptive loop. More recently, also recovery-based methods for polynomial splines over T-meshes have been proposed \citep{Anitescu2018}. 
In the context of (truncated) hierarchical B-splines, we mention \citep{Buffa2016,Gantner2017} where residual-based element-wise estimators have been analyzed for linear second-order elliptic problems. These concepts have been extended in \citep{Buffa2018a} for function-based estimators.

The main goal of this paper is to introduce, building upon the work of \citep{Bank1993,Vuong2011} on multi-level estimators, a simple \textit{a posteriori} error estimator for linear fourth-order elliptic partial differential  equations (PDEs) of engineering relevance, such as Kirchhoff plates and Kirchhoff-Love shells, by making use of hierarchical B-splines. The proposed method is based on the solution of an additional residual-like problem which is formulated on a so-called \textit{bubble space}. Given the locality of the functions in the aforementioned space, the resulting linear system is in general small and easy to solve. We remark that by exploiting this technique we bypass completely the computational burden stemming from the evaluation of higher order derivatives and their jumps across edges, which would be required in classical residual-based error indicators. This is particular advantageous for Kirchhoff-Love shells, where the terms associated to higher derivatives involve the cumbersome and error-prone task of taking (nested) covariant derivatives of quantities of interest, see for instance \citep{Maurin2018}. Indeed, the use of residual-based error estimators for Kirchhoff-Love shells is computationally expensive and an efficient implementation is far from being trivial. Moreover, contrary to residual-based techniques, we do not require the computation of additional integrals which take into account the jumps of derivatives across element interfaces. These features play again in favor of the computational efficiency and simplicity of the proposed methodology. 

A known drawback of this family of estimators is that their reliability and efficiency is subjected to the saturation assumption onto the underlying augmented space. As noted in \citep{Bank1993}, this assumption is, in general, problem dependent and it can, potentially, become crucial. However, we highlight that our focus is to find a good indicator for steering adaptive simulations in the scope of structural mechanics and all our numerical experiments confirm that the proposed estimator performs well in this task. Furthermore, it is worth noting that in all our tests the estimator provides also an excellent approximation of the true error. 

The structure of the paper is as follows. \Cref{sec:IGA} introduces the fundamental concepts needed for hierarchical B-splines. \Cref{sec:formulations} briefly recalls the Kirchhoff plate and Kirchhoff-Love formulations. Then, \Cref{sec:bubbles} discusses the proposed bubble error estimator and a possible algorithm. \Cref{sec:numericalExamples} presents several numerical examples. Finally, some conclusions are drawn in \Cref{sec:conclusions}.
\newcommand\superHB{\scriptscriptstyle\mathcal{HB};\;}
\newcommand\superTHB{\scriptscriptstyle\mathcal{THB};\;}
\newcommand{\Norm}[1]{\left\lVert#1\right\rVert_2}
\newcommand{\delom}{{\;\text{d}\Omega}}
\newcommand{\delg}{{\;\text{d}\Gamma}}
\newcommand{\energyNorm}[1]{\left\lVert#1\right\rVert_{E(\Omega)}}
\newcommand{\snot}{{S}}
\newcommand{\s}{{\snot \left( \T, p \right) }}
\newcommand{\T}{{\mathcal{T}}}
\newcommand{\op}{{\Omega\textsuperscript{\text{phy}}}}
\newcommand{\of}{{\Omega\textsuperscript{\text{fict}}}}


\section{IsoGeometric Analysis} \label{sec:IGA}

\newcommand{\bg}{}
In this Section, the notation and basic concepts and definitions related to B-splines and NURBS are reviewed, following closely the derivation in~\citep{Dangella2018}. For further details, the reader is referred to~\citep{Piegl1995, Cottrell2009}, and references therein.

A B-spline basis function of degree $p$ is generated starting from a non-decreasing sequence of real values that corresponds to a set of coordinates in the parameter space referred to as knot vector, denoted in the following as $\Xi$.
Given a knot vector, a univariate B-spline basis function $b_{i,p}$ can be defined recursively using the Cox-de Boor algorithm~\citep{Piegl1995}.
It is worth mentioning that the smoothness of the obtained B-spline basis is $C^{p-k}$ at every knot, where $k$ denotes the multiplicity of the considered knot, while it is $C^\infty$ everywhere else.

%
The definition of multivariate B-splines is achieved in a straight-forward manner using the tensor product of univariate B-splines. Indeed, the multivariate B-splines $\mathcal{B}_{\mathbf{i},\mathbf{p}}$ are obtained as:
\begin{equation*}
\mathcal{B}_{\mathbf{i},\mathbf{p}}=\prod_{r=1}^{d_{r}}b_{i_{r},p_{r}}^{r}
\end{equation*}
with $r=1,...,d_{r}$, where $d_{r}$ is the dimension of the parameter space and $p_{r}$ denotes the polynomial degree in the parametric
direction $r$, respectively. Additionally, the multi-index $\mathbf{i}=\left\{ i_{1},...,i_{d_{r}}\right\} $
denotes the position in the tensor product structure and $\mathbf{p}=\left\{ p_{1},...,p_{r}\right\} $
indicates the vector of polynomial degrees. 

Lastly, let us define a B-spline parametrization $\mathbf{S}$ as a linear combination of multivariate B-spline basis functions and control points as follows:
\begin{equation*}
\mathbf{S}=\sum_{\mathbf{i}}\mathcal{B}_{\mathbf{i},\mathbf{p}}\mathbf{P}_{\mathbf{i}} \, ,
\end{equation*}
where the coefficients $\mathbf{P}_{\mathbf{i}}\in\mathbb{R}^{d}$ of the linear combination are the control points and $d$ denotes the dimension of the physical space. We remark that all the concepts summarized here can be readily transferred to NURBS, for further details we refer to~\citep{Cottrell2009}.
In the rest of the paper, without loss of generality, the degree $ p $ will be considered equal in each parametric direction and will be omitted from the notation.

\subsection{Hierarchical B-splines} \label{subsection:hierchical}

In this Subsection, the concept of hierarchical B-spline basis, denoted by $ \mathcal{HB} $, is introduced. This allows us to build a basis that is locally refinable and therefore to overcome the limitations intrinsic to the tensor-product nature of B-splines and NURBS. 

Let $ V^0 \subset V^1 \subset \cdots \subset V^{N} $ be a sequence of nested spaces of splines defined on a parametric domain $ \hat{\Omega}$, where, to keep the notation simple, we consider the one-dimensional case. Each space $ V^l $, $ l=0 ,\ldots, N $ is spanned by the B-spline basis $ \mathcal{B}^l $ of degree $ p $, associated to level $ l $ and the corresponding knot vector $\Xi^l$. 

Our task now is to identify a set of analysis-suitable functions $ \mathcal{N} \subset \bigcup_{l} \mathcal{B}^l $. To this end, let us define as isogeometric elements $Q$ a partition of $ \hat{\Omega} $. In particular, given $\hat n^l \in \mathbb{N}$, let $ \hat \Xi^l = (\hat\xi^l_0, \ldots, \hat \xi^l_{\hat n^l} ) $ be the knot vector composed of non-decreasing knots of $ \Xi^l $ without repetition, and let:
\begin{align*}
\mathcal{Q}^l = \left\lbrace \mathcal{Q}^l_i \;  | \; \mathcal{Q}^l_i = ( \hat \xi^l_i, \hat \xi^l_{i+1} ) , i=0,\ldots, \hat n^l-1 \right\rbrace 
\end{align*}
be the set of open intervals constituting the non-empty knot spans of $ \Xi^l $. The elements of the multi-level mesh can be any non-overlapping partition $ Q \subset \bigcup_{l} \mathcal{Q}^l $ of $ \hat{\Omega} $, such that $\bigcup_{\epsilon \in Q} \bar \epsilon = \hat{\Omega}$.
%
%
Moreover, let $ Q^l = Q \cap \mathcal{Q}^l $ be the elements of level $ l $, and let $ \hat{\Omega}^l = \bigcup_{\epsilon \in Q^l} \bar \epsilon $ denote their domain. Furthermore, let $ \hat{\Omega}^l_+ = \bigcup_{l^*=l+1}^N  \hat{\Omega}^{l^*} $ be the refined domain with respect to level $ l $ and, analogously, let $ \hat{\Omega}^l_- = \bigcup_{l^*=0}^{l-1}  \hat{\Omega}^{l^*} $ be the coarser domain. 
%
\noindent Given a set of elements $Q$, we still need to define a set of basis functions $ \mathcal{N} \subset \bigcup_{l} \mathcal{B}^l $ suitable for analysis. To this end, we consider the set of functions with support on the elements of level $l$ as the set of \emph{active} functions $ \mathcal{B}^l_a =  \left\lbrace b \; |  b \in \mathcal{B}^l, \, \supp{b} \cap \hat{\Omega}^l \neq \emptyset \right\rbrace \subset\mathcal{B}^l $. Among these, a subset of linearly independent functions has to be chosen. Following~\citep{Hennig2016,Dangella2018}, we partition $ \mathcal{B}^l_a $ into:
\begin{align*}	
	\mathcal{B}^l_- &=  \left\lbrace b \; |\;  b \in \mathcal{B}^l_a, \, \supp{b} \cap \hat{\Omega}^l_- \neq \emptyset \right\rbrace \, , \\
	 \mathcal{B}^l_+ &=  \left\lbrace b \; |\;  b \in \mathcal{B}^l_a, \, \supp{b} \cap \hat{\Omega}^l_+ \neq \emptyset \right\rbrace \setminus \mathcal{B}^l_- \, , \\
 \mathcal{B}^l_= &=  \left\lbrace b \; |\;  b \in \mathcal{B}^l_a, \, \supp{b}  \subset \hat{\Omega}^l  \right\rbrace  \, .
\end{align*}
%
%
Then, the hierarchical B-spline basis $ \mathcal{HB} $~\citep{Vuong2011,Kraft1997} is defined as follows:
  \begin{align}
	\label{eq:hb}
	\mathcal{HB} &= \bigcup\limits_l \; \mathcal{HB}^{l}\, ,  \\
	\mathcal{HB}^{l} &= \left( \mathcal{B}^l_= \cup \mathcal{B}^l_+ \right). \nonumber
\end{align}
Namely, $ \mathcal{HB} $ is the set of B-splines of each level $ l $ whose support intersects only elements of level $ \tilde{l} \geq l $ and at least one element of level $ l $. It was proven in~\citep{Vuong2011} that this set is composed of linearly independent functions and therefore suitable for the analysis.


\section{Formulations} \label{sec:formulations}

In this section, the Kirchhoff plate and the Kirchhoff-Love shell formulations are stated in the context of IGA. The reader is referred to~\citep{Reali2015,Kiendl2009,Cirak2006,Maurin2018} for a more comprehensive review of the formulations.

\subsection{The Kirchhoff plate formulation} \label{subsec:plateFormulation}

Here, following the derivation in~\citep{Reali2015}, we introduce the Kirchhoff plate problem, governed by the bi-Laplace differential operator. Let us define an open set $\Omega \subset \mathbb{R}^2$ with a sufficiently smooth boundary $\partial \Omega$, such that the normalized normal vector $\boldsymbol{d}$ is well-defined (almost) everywhere.   
Additionally, we assume that the boundary $\Gamma = \partial \Omega$ can be decomposed into $\Gamma =\overbar{\Gamma_w \cup \Gamma_Q}$ and $\Gamma =\overbar{\Gamma_{\phi} \cup \Gamma_{M}}$, such that $\Gamma_w \cap \Gamma_Q = \varnothing$ and $\Gamma_{\phi} \cap \Gamma_{M} = \varnothing$, respectively. We formulate the strong form of the problem as follows:
\begin{alignat}{2}\label{eq:strongBVPPlate}
	D \varDelta^2 u &= g \quad &&\text{in} \quad \Omega \\ 
	u &= u_{\space \Gamma} \quad &&\text{on} \quad \Gamma_w \nonumber \\ 
	- \nabla u \cdot \boldsymbol{d} &= \phi_{\space \Gamma} \quad &&\text{on} \quad \Gamma_{\phi} \nonumber \quad ,\\
	\nu D \varDelta u + (1 - \nu) D \, \boldsymbol{d} \cdot (\nabla \nabla u)\boldsymbol{d} &= M_{\space \Gamma} \quad &&\text{on} \quad \Gamma_{M} \nonumber \\
	D ( \nabla(\varDelta u) + (1 - \nu) \, \boldsymbol{\varPsi}(u)\,) \cdot \boldsymbol{d} &= Q_{\space \Gamma} \quad &&\text{on} \quad \Gamma_{Q}  \nonumber
\end{alignat} 
where $u$ represents the deflection of the plate, $D$ its bending stiffness, $\nu$ is the Poisson ratio, $g$ is the load per unit area in the thickness direction, $u_{\space \Gamma}$, $\phi_{\space \Gamma}$, $M_{\space \Gamma}$ and $Q_{\space \Gamma}$ are the prescribed deflection, rotation, bending moments and effective shear, respectively. The bending stiffness $D$ is defined as:
\begin{align*}
D = \frac{E t^3}{12(1 - \nu^2)} \, ,
\end{align*}
where $E$ is the Young modulus and $t$ denotes the thickness of the plate, which without loss of generality we suppose to be a constant in $\Omega$.
Using variational calculus, the weak form of problem \eqref{eq:strongBVPPlate} reads: find $u  \in \: V$ such that
 \begin{align}\label{eq:modelWeakFormPlate}
 &a(u,v) = F(v)  \qquad \forall v \in V \, ,
 \end{align}
where $V \subset H^2(\Omega)$ denotes an infinite-dimensional space that depends in general on the boundary conditions of the problem at hand, for further details we refer to \citep{Ciarlet2002}. 
For simplicity and without loss of generality, in the following we assume the Poisson ratio $\nu$ to be zero and we choose $E$ and $t$ such that $D = 1$. Consequently, the bilinear form $a: V \times V \rightarrow \mathbb{R}$ can be expanded as follows:
\begin{align*}
 	a(u,v) &= \int_{\Omega} \nabla (\nabla v) \, : \, \nabla (\nabla u) \, \text{d}\Omega \, ,
\end{align*}
and similarly the linear functional $F: V \rightarrow \mathbb{R}$ reads:
 \begin{align*}
 	F(v) &= \int_{\Omega} g v \, \text{d}\Omega \, .
\end{align*}
It can be shown that \eqref{eq:modelWeakFormPlate} satisfies the conditions of the Lax-Milgram theorem, therefore $u \in V$ is the unique solution of the variational problem \eqref{eq:modelWeakFormPlate}. 

Finally, we discretize the weak form \eqref{eq:modelWeakFormPlate} employing the \textit{Bubnov-Galerkin} method and using as basis the hierarchical B-splines $\mathcal{HB}$. Let $V_h \subset V$ be a finite dimensional subspace defined as:
\begin{align*}
V_h = \text{span }\lbrace b \circ \mathbf{S}^{-1} \vert \, b \in \mathcal{HB} \rbrace.
\end{align*}
where $\mathbf{S}$ is the spline parametrization defined in~\Cref{sec:IGA}. We remark that due to the requirements on the discrete admissible space for the deflection, the basis functions must be at least $C^1$ continuous globally such that the bilinear form is well-defined. This requirement is easily fulfilled within one patch by the use of (hierarchical) B-splines of maximum continuity (of degree $p \geq 2$). Then, the discrete weak form can be written as: find $u_h  \in \: V_h$ such that 
\begin{align}\label{eq:modelWeakFormPlate_discrete}
 &a(u_h,v_h) = F(v_h)  \qquad \forall v_h \in V_h \, ,
\end{align}
where $u_h$ represents the discrete solution defined as:
\begin{align*}
u_h = \sum_{b \, \in \, \mathcal{HB}} b \, c_b \, ,
\end{align*}  
with $c_b \in \mathbb{R}$ denoting the unknown control variables and $b \in \mathcal{HB}$ the corresponding hierarchical basis introduced in~\Cref{sec:IGA}. 
 
\subsection{The Kirchhoff-Love shell formulation}


In the following, the Kirchhoff-Love formulation in its weak form is briefly introduced, following closely the derivation and the notation used in~\citep{Kiendl2009,Cirak2006}.

\subsubsection{The weak formulation}

Similarly to what has been described in~\Cref{subsec:plateFormulation} for the Kirchhoff plate, starting from the strong formulation of the Kirchhoff-Love problem (for further details see~\citep{Maurin2018}) and applying variational calculus, the weak form reads: find $\boldsymbol{u}  \in \: V$ such that
\begin{align}\label{eq:modelWeakForm}
	&a(\boldsymbol{u},\boldsymbol{v}) = F(\boldsymbol{v})  \qquad \forall \boldsymbol{v} \in V \, ,
\end{align}
where the choice of space $V \subset H^2(\Omega)$ depends on the boundary conditions of the problem at hand.
Moreover, $a:V \times V \rightarrow \mathbb{R}$ is a continuous and strongly coercive bilinear form and $F: V \rightarrow \mathbb{R}$ is a continuous linear functional. They can be expanded, respectively, as follows:
\begin{subequations}\label{eq:membersWeakForm}
\begin{align*}
	a(\boldsymbol{u},\boldsymbol{v}) &= \int_{\Omega} \boldsymbol{\varepsilon}(\boldsymbol{v}) \, : \, \boldsymbol{n}(\boldsymbol{u}) \delom + \int_{\Omega} \boldsymbol{\kappa}(\boldsymbol{v}) \, : \, \boldsymbol{m}(\boldsymbol{u}) \delom \, , \\
	F(\boldsymbol{v}) &= \int_{\Omega} \boldsymbol{v} \cdot \boldsymbol{b} \delom + \int_{\Gamma_N} \boldsymbol{v} \cdot \boldsymbol{p} + \boldsymbol{\omega}(\boldsymbol{v}) \cdot \boldsymbol{r} \delg \, ,
\end{align*}
\end{subequations}
where $\boldsymbol{\varepsilon}$, $\boldsymbol{\kappa}$ denote the membrane and bending strain tensors, respectively, and $\boldsymbol{n}$, $\boldsymbol{m}$ are their energetically conjugate stress resultants. Additionally, $\boldsymbol{b}$ is the applied body load, $\boldsymbol{p}$ and $\boldsymbol{r}$ represent the applied traction force and traction moment, respectively. For a detailed derivation, we refer to~\citep{Kiendl2009,Apostolatos2015}.
Note that the weak form \eqref{eq:modelWeakForm} satisfies the conditions of the Lax-Milgram theorem and therefore the solution $\boldsymbol{u} \in V$ exists and it is unique.

Employing the IGA paradigm, we utilize the same (hierarchical) B-splines and NURBS basis functions for the geometry description and for the field approximation. 
Let us recall the definition of the discrete space:
\begin{align*}
V_h &= \text{span }\lbrace b \circ \mathbf{S}^{-1} \vert \, b \in \mathcal{HB} \rbrace \, ,
\end{align*}
by employing the \textit{Bubnov-Galerkin} discretization of \eqref{eq:modelWeakForm}, we get the following discrete weak formulation of the Kirchhoff-Love shell problem: find $\boldsymbol{u}_h \in V_h$ such that
\begin{align} %
\label{eq:modelFEForm}
& a\left(\boldsymbol{u}_h ,\boldsymbol{v}_h \right) = F \left(\boldsymbol{v}_h \right)  \qquad \forall \boldsymbol{v}_h \in V_h \, .
\end{align}
Again, it is worth noting that due to the requirements on the discrete admissible space for the displacement field, the basis functions must be at least $C^1$ continuous globally such that the bending operator is well-defined. 

Finally, we remark that in most cases, due to efficiency, the parametrization $\mathbf{S}$ can be expressed in terms of tensor product B-splines of the coarsest level $\mathcal{B}^0$ instead of the hierarchical basis $\mathcal{HB}$ \citep{Garau2018}.

\newcommand{\picsDir}{pictures/numericalExamples/pics}
\newtheorem{remark}{Remark}

\section{A posteriori error estimator for Kirchhoff plates and Kirchhoff-Love shells} \label{sec:bubbles}

In the following Section we introduce a variant of the error estimator analyzed in \citep{Bank1993} and we extend its isogeometric version, proposed in \citep{Juettler2010}, to fourth-order PDEs of mechanical relevance, such as Kirchhoff plates and Kirchhoff-Love shells. We choose these models since, due to the Kirchhoff hypotheses, we avoid a priori the negative effects related to shear locking. Then, later in the Section, we show a possible implementation of the proposed indicator in the context of hierarchical B-splines.

\subsection{The bubble error estimator}
\begin{figure}[!h]
	\begin{subfigure}[t]{0.475\textwidth}
	\centering
	\includegraphics[scale=1.0]{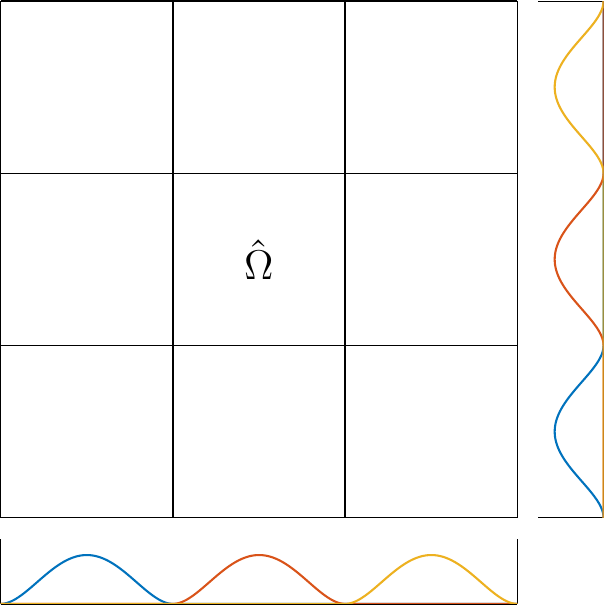}
	\caption{Bubble basis $p=4$.}
	\label{fig:bubble_basis_p4}
	\end{subfigure}
	\begin{subfigure}[t]{0.475\textwidth}
	\centering
	\includegraphics[scale=1.0]{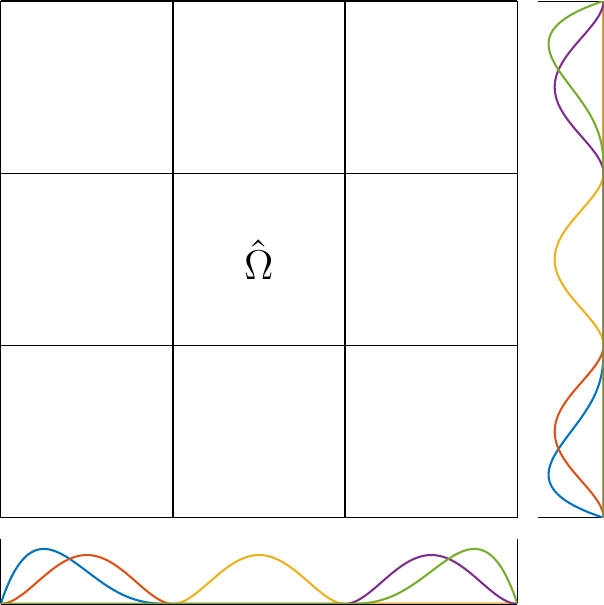}	
	\caption{Bubble basis $p=4$ with boundary functions.}
	\label{fig:bubble_basis_p4_withBCs}
	\end{subfigure}		
\par\bigskip
	\begin{subfigure}[t]{0.475\textwidth}
	\centering
	\includegraphics[scale=1.0]{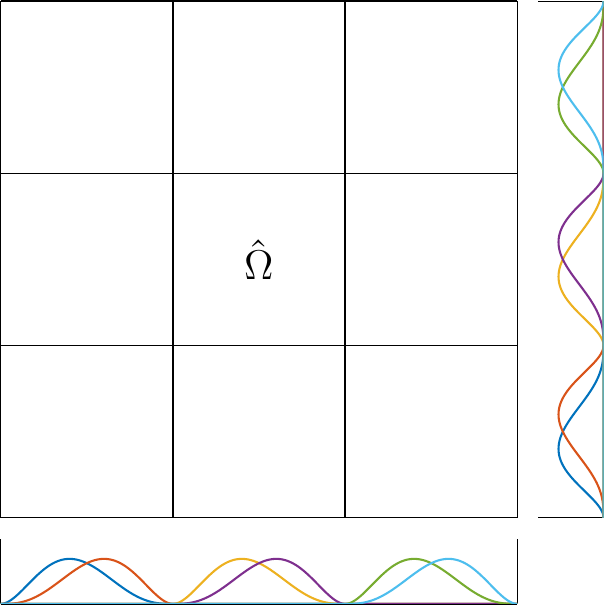}	
	\caption{Bubble basis $p=5$.}
	\label{fig:bubble_basis_p5}
	\end{subfigure}
	\begin{subfigure}[t]{0.475\textwidth}
	\centering
	\includegraphics[scale=1.0]{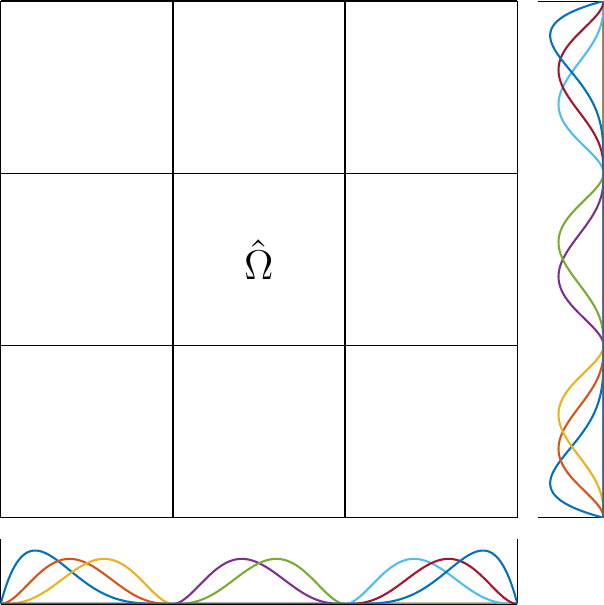}		
	\caption{Bubble basis $p=5$ with boundary functions.}
	\label{fig:bubble_basis_p5_withBCs}
	\end{subfigure}	
	\caption{Example of construction of bubble functions $b_h$ as tensor product of univariate functions for the case $p=3,4$, where on the left-side only internal functions are depicted whereas on the right-side also boundary functions are plotted.}
	\label{fig:bubble_basis}
\end{figure}

This family of \textit{a posteriori} estimators was introduced in \citep{Bank1993} in the context of the p-version of the Finite Element Method and has been successfully applied to T-splines in \citep{Juettler2010} for second-order boundary value problems. To introduce the concept in the context of hierarchical IGA, let us define the finite dimensional solution space $V^p_h$ as the span of hierarchical B-splines basis functions of order $p$, where we restrict our analysis to the case $p > 2$ (see Remark~\ref{rem:2} for the case $p =2$). Then, let us recall the discrete solution $u_h \in V^p_h$, where without loss of generality we make no distinction between the plate or the shell problem. We now assume there exists a larger space $V^p_h \subset \widetilde{V}^p_h \subset V$ such that the following decomposition holds:
\begin{align*}
\widetilde{V}^p_h = V^p_h \oplus W^{p+1}_h \, ,
\end{align*} 
where $W^{p+1}_h$ is the space in which we seek a good estimation of the error $e_h \approx e = \left\lVert u - u_h \right\rVert$, in a suitable norm $\left\lVert \cdot \right\rVert$. In particular, recalling some definitions provided in Section \ref{sec:IGA}, let us characterize $W^{p+1}_h$ as follows:
\begin{align*}
W^{p+1}_{h} = \bigcup_{l=0}^N W^{p+1}_{h,l} \, ,
\end{align*}
where 
\begin{align*}
W^{p+1}_{h,l} = \text{span }\lbrace b \circ \mathbf{S}^{-1} \vert \, b \in \mathcal{B}^l_{a,p+1} \cap C^1 ( \overbar{\hat{\Omega}^l} )  \text{ and } \exists \, \epsilon \subseteq \hat{\Omega}^l : \text{ supp}(b) = \overbar{\epsilon} \rbrace.
\end{align*}
In words, $W^{p+1}_{h,l}$ is the space spanned by active B-splines (or better in this case, Bernstein polynomials) of level $l$ and degree $p+1$ such that their support is compact and overlaps exactly with one active element $\epsilon$ of level $l$, e.g. see~\Cref{fig:bubble_basis} where one level is depicted.
We are now ready to define the a posteriori error estimate as: find $e_h  \in \: W^{p+1}_h$ such that
 \begin{equation}\label{eq:bubbleSystem}
 a(e_h,b_h) = F(b_h) - a(u_h,b_h)  \qquad \forall b_h \in W^{p+1}_h  \, ,
 \end{equation}  
where $b_h$ is referred to in the following as \textit{bubble functions}. We remark that we require these functions to be at least $C^1$ on their support such that the corresponding bilinear form $a$ is well-defined. Notice that due to the compact support of $b_h$, the linear system corresponding to the discrete error weak form \eqref{eq:bubbleSystem} is block diagonal and in general easy to solve. 
Moreover, we highlight that the evaluation of the bilinear form $a$ appearing in~\eqref{eq:bubbleSystem} follows the same steps as the original plate/shell problem, where different test functions $b_h$ are used. Hence, setting up the bubble problem requires only minor changes into an existing isogeometric code.    

Once we have found the unknown coefficients $e_h$, we define the element-wise error estimator $\eta_{\epsilon}$ as follows:
\begin{equation}\label{eq:bubbleSystem_bis}
\eta_{\epsilon} = C_a \left\lVert e_h \right\rVert_{E(\epsilon)}  \quad \forall \epsilon \in Q \, ,
\end{equation}  
where we recall that $\epsilon$ is an isogeometric element of the hierarchical mesh $Q$ and $\left\lVert \cdot \right\rVert_{E(\epsilon)}$ denotes the energy norm restricted to element $\epsilon$. Additionally, we introduce $C_a$ as an empirical constant which, considering all our numerical experiments, we claim is independent from the chosen degree and from the problem at hand. In the following, we set $C_a = 3$ once and for all in order to shift the estimator above the true error.

At this point, it is worth highlighting that the proposed error estimator is computationally inexpensive, easy to implement and is embarrassingly parallelizable due to the choice of the disjoint bubble space $W^{p+1}_h$. A possible algorithm to estimate the error is summarized in Algorithm~\ref{alg:estimate}, where we design the code to solve the additional problem \eqref{eq:bubbleSystem} level-wise over the hierarchical mesh $Q$.

\begin{algorithm}[H] 
\begin{algorithmic}[1]
	\Procedure{Estimate error }{numerical solution $u_h$, hierarchical mesh $Q$} 
	\State Initialize vector $\eta$
	\For{\textbf{each} level $l$ of $Q$}
		\State Build bubble space $W^{p+1}_{h,l}$
		\State Given $u_h$, get $u_h^l$ of level $l$
		\State Given $u_h^l$, set up and solve the additional system \eqref{eq:bubbleSystem}
		\For{\textbf{each} $\epsilon \in Q^l$}
			\State Compute the element-wise indicator $\eta_{\epsilon}$ \eqref{eq:bubbleSystem_bis}
			\State Store $\eta_{\epsilon}$ in $\eta$
		\EndFor
	\EndFor	
	\State Return $\eta$
	\EndProcedure
\end{algorithmic} 
\caption{Bubble error estimator algorithm}\label{alg:estimate}
\end{algorithm}

\begin{remark}[Inhomogeneous boundary conditions]
On one hand, we highlight that whenever inhomogeneous boundary conditions of Neumann-type are applied to the problem at hand, problem \eqref{eq:bubbleSystem} must also contain the corresponding additional terms in the right-hand-side. In these cases, we augment the bubble space $W^{p+1}_h$ with suitable boundary bubbles to properly capture the error on the imposition of natural boundary conditions. An example is depicted in~\Cref{fig:bubble_basis} for bubbles of degree $p=4,5$. On the other hand, we make the assumption that the error on the imposition of inhomogeneous Dirichlet-type boundary conditions is negligible and therefore no additional shape functions are needed for these cases.
\end{remark}

\begin{remark}[The case $p=2$] \label{rem:2}
We remark that isogeometric Kirchhoff-Love elements notoriously suffer from membrane locking phenomena, specially for low degree discretization \citep{Kiendl2009}, and it can be considered common practice to use the coarsest (and lowest degree) description of the geometry that is still exact while employing a \textit{k}-refined space for the solution field to mitigate these effects. 
Nevertheless, for the sake of completeness, we should highlight that the construction of the bubble space $W^{p+1}_h$ is not straightforward in the case where the underlying B-splines are quadratic. 
 \noindent We propose as possible remedies:
\begin{enumerate}
\item Construction of an additional sub-grid obtained as dyadic refinement of the original hierarchical mesh and definition of the bubble functions on the corresponding space defined by the new grid.

\item Introduction of suitable edge functions in the bubble space.
\end{enumerate}
Notice that, although both approaches could be feasible, they add a non-negligible complexity to the method. In the first case, a new, finer hierarchical grid has to be created at every iteration of the adaptive loop and numerical integration of \eqref{eq:bubbleSystem} has to be performed over it. In the second case, the definition of which edge bubbles have to be considered in the hierarchical case is not unique and considerably worsen the simplicity of the method. Indeed, the corresponding linear system would become more difficult to solve since the element-wise locality of the bubble functions would be lost. Therefore, 
in light of these considerations, we will not address further the case $p=2$ in the scope of this work.     
\end{remark}

\begin{remark}[Jumps in the residual]
In the last remark, we notice that, although the residual exhibits jumps across element edges for the discretization of degree $p=3$, in all our computations the estimator behaves still optimally even if the corresponding jump terms are not accounted for in \eqref{eq:bubbleSystem}.  
\end{remark}

\subsection{Marking and refine strategies}

Here, we briefly recall that once an estimate of the error element-wise is available, there exist several strategies proposed in the literature to mark elements for refinement. In all our numerical examples we use the so-called \textit{maximum strategy}, which can be summarized as follows. Let $\gamma \in (0,1)$ be a user-defined threshold, all elements such that
\begin{align*}
\eta_{\epsilon} > \gamma \, \tilde{\eta}_{\varepsilon} \, , \quad \text{where} \quad \tilde{\eta}_{\varepsilon} = \max_{\varepsilon \in Q} \, \eta_{\varepsilon} \, ,
\end{align*}
are marked for refinement. In the following, if not stated otherwise, we will always use $\gamma = 0.5$.
Additionally, since for basis functions of order $p$ their support usually spans up to $p+1$ elements, we try to reduce the number of isolated marked elements. Indeed, in general, isolated refined elements do not improve the solution space $V^p_h$ (they do not add any additional degree of freedom) but only the associated numerical integration is refined. To avoid that, we slightly modify the maximum strategy marking algorithm such that once an element $\epsilon$ of level $l$ is marked for refinement, the algorithm also marks all the neighbors of $\epsilon \in Q^l$, where $Q^l$ is the set of active elements of level $l$.
Finally, we remark that our marking algorithm is also designed to preserve the admissibility (as defined in~\citep{Buffa2016a}) of the hierarchical mesh between consecutive iterations of the adaptive loop. We recall that admissible meshes are guaranteed to have a bounded number of basis functions acting on any element of the mesh, such that interaction between very fine and very coarse functions is avoided. 
If not stated otherwise, we will set the class of admissibility $m$ to be $m = p-1$ in all our numerical examples. 
For the sake of brevity we do not discuss the details here, the reader is referred to \citep{Buffa2016a,Buffa2017,Bracco2018} for a comprehensive review of the concept of admissibility.

\newcommand{\graphDir}{pictures/numericalExamples/graphs}
\newcommand{\dataDir}{pictures/numericalExamples/data}

\section{Numerical Examples}
\label{sec:numericalExamples}

All the numerical experiments presented in the following section have been implemented in the open-source and free Octave/Matlab package \textit{GeoPDEs} \citep{Vazquez2016}, which is a software suite for the solution of partial differential equations specifically designed for isogeometric analysis.  	

\subsection{Kirchhoff plate}

In this Section we present several adaptive computations in the context of Kirchhoff plates, thought to give a first assessment of the performance of the proposed error estimator. We demonstrate the applicability of our method both for problems which exhibit smooth and singular solutions. In all cases presented here, the proposed error estimator shows excellent performance in steering the adaptive simulation, yielding the expected optimal rates of convergence in the asymptotic regime. 

\subsubsection{Smooth solution on a square plate}

\begin{figure}
	\begin{subfigure}[t]{0.475\textwidth}
		\centering
	\includegraphics[scale=1.0]{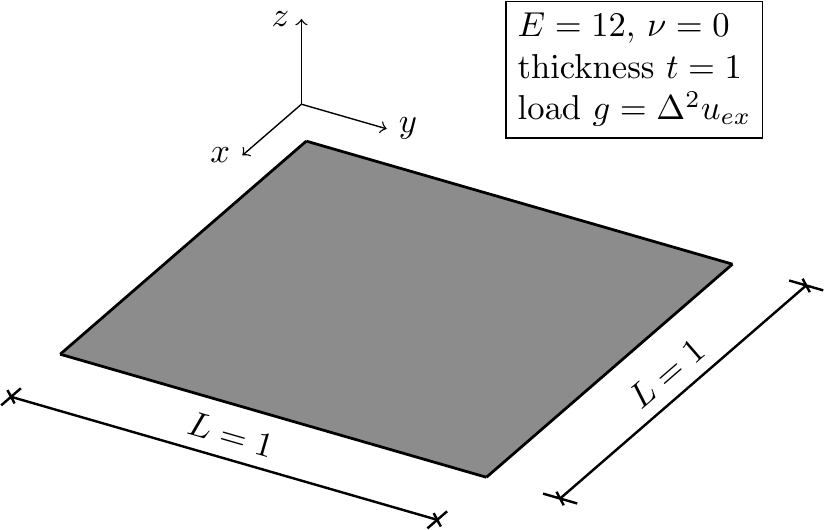}
	\caption{Geometry description of the problem.}
	\end{subfigure}
	\hfill
	\begin{subfigure}[t]{0.475\textwidth}
		\centering
	\includegraphics[scale=1.0]{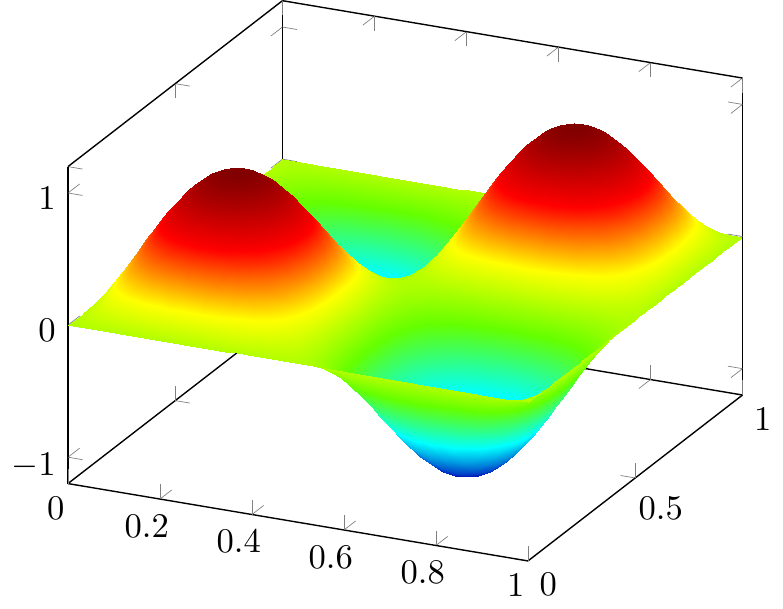}
	\caption{Exact solution $u_{ex} = \sin(2 \pi x) \sin(2 \pi y)$.}
	\end{subfigure}
	\caption{Geometry, physical parameters and exact solution contour of the square plate example.}
	\label{fig:plate_geo}
\end{figure}

In the first example, we analyzed the behavior on the bubble error estimator compared to a classical residual-based type error estimator. We define the computational domain to be the unit square $\Omega = [0,1]^2$ and we impose the following homogeneous boundary conditions:
\begin{alignat*}{2}
u &= 0 \qquad &&\text{on} \quad \partial \Omega  \quad .\\
\nu D \varDelta u + (1 - \nu) D \, \boldsymbol{d} \cdot (\nabla \nabla u)\boldsymbol{d} &= 0 \qquad &&\text{on} \quad \partial \Omega 
\end{alignat*}  
Additionally, the applied load $g$ is constructed such that it fulfills the following manufactured solution $u_{ex} = \sin(2 \pi x) \sin(2 \pi y)$. Namely, $g$ reads as follows:
\begin{align*}
g = 64 \pi^4 \sin(2 \pi x) \sin(2 \pi y). 
\end{align*}
The structure, physical parameters and the exact solution are depicted in Figure~\ref{fig:plate_geo}.
Notice that, as predicted by classical \textit{a priori} error estimates (see for instance \citep{Ciarlet2002,Quarteroni2008} or in the context of IGA \citep{Bazilevs2006}) and since $u_{ex}$ is sufficiently regular, we expect the following convergence rate:
\begin{align}
\left\lVert u_{ex} - u_h \right\rVert_{H^s(\Omega)} \leq C h^{p+1-s} \left\lVert u_{ex} \right\rVert_{H^{p+1}(\Omega)} \qquad p+1 > s \, ,
\label{eq:apriori_smooth}
\end{align}
where $p$ denotes the polynomial degree of the hierarchical basis and $h$ represents the maximum mesh size, namely:
\begin{align*}
h = \max_{\epsilon \, \in \, Q} h_\epsilon.
\end{align*}
Additionally, in all our numerical examples, we choose $s=2$. Moreover, for simplicity and without loss of generality, we measure the error in the $H^2$ semi-norm, which corresponds to the energy norm of the problem at hand. It is worth noting that in most cases we will use the number of degrees-of-freedom (dofs) instead of the mesh size $h$, since we can relate these two quantities via $h \approx \text{dofs}^{- 1/d}$, where $d$ denotes the dimensionality of the problem.

\begin{figure}
	\centering
	\begin{subfigure}[t]{0.495\textwidth}
		\centering
	\input{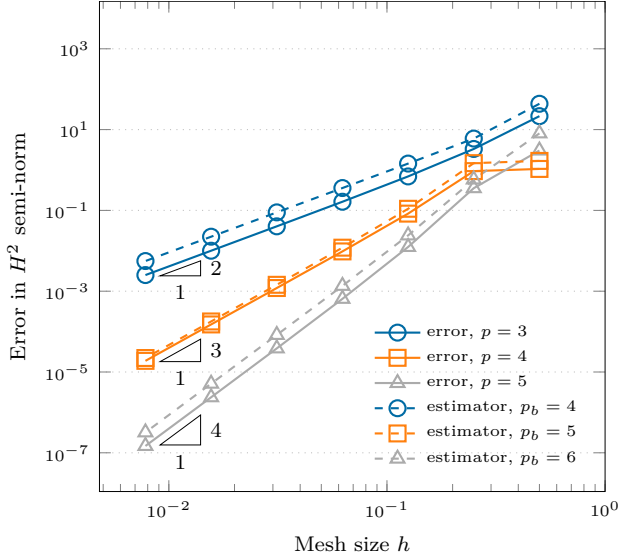}
		\caption{Bubble estimator.}
	\end{subfigure}
	\hfill
	\begin{subfigure}[t]{0.495\textwidth}
		\centering
	\input{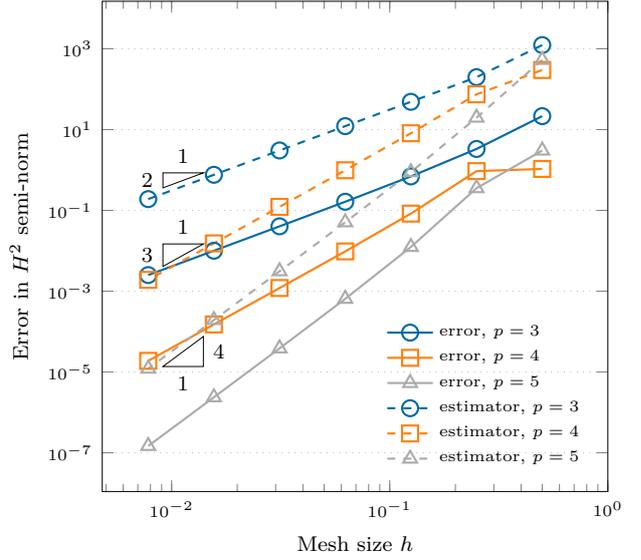}
		\caption{Residual-based estimator.}
	\end{subfigure}
	\caption{Study of the convergence of the error measured in the $H^2$ semi-norm, bubble and residual-based estimators against the mesh size $h$ for different $p$ on a smooth plate problem ($u_{ex}=\sin(2 \pi x) \sin(2 \pi y)$). Uniform refinement.}
\label{fig:convergence_smooth_GR}
\end{figure}

For the sake of comparison, in the case of Kirchhoff plates only, we drive the adaptive method by means of both the residual-based and the bubble-based estimators.
To test our implementation, we perform at first uniform refinement for different degrees $p=3,4,5$ and check the convergence rate of the error in the energy norm and both estimators against the element size $h$, as depicted in~\Cref{fig:convergence_smooth_GR} for the bubble and residual-based, respectively. In all the presented cases an optimal asymptotic rate of convergence is observed, both for the error and the estimator. Furthermore, we remark that the bubble estimator is much closer to the real error compared to the residual-based.
%

\begin{figure}
	\centering
	\begin{subfigure}[t]{0.495\textwidth}
		\centering
	\input{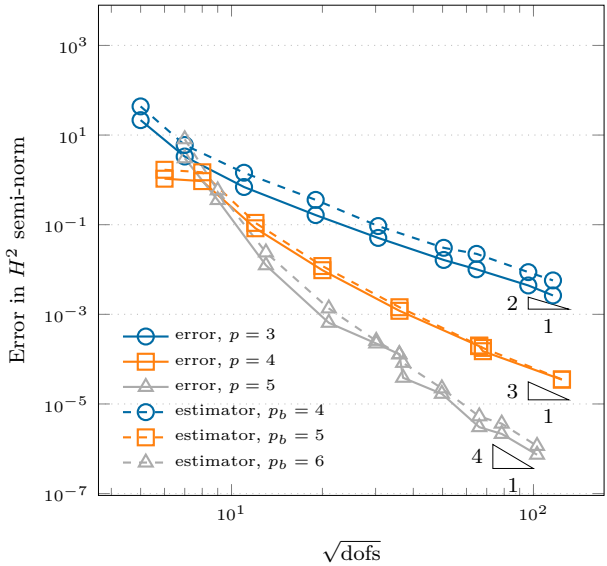}
		\caption{Bubble estimator.}
	\end{subfigure}
	\hfill
	\begin{subfigure}[t]{0.495\textwidth}
		\centering
	\input{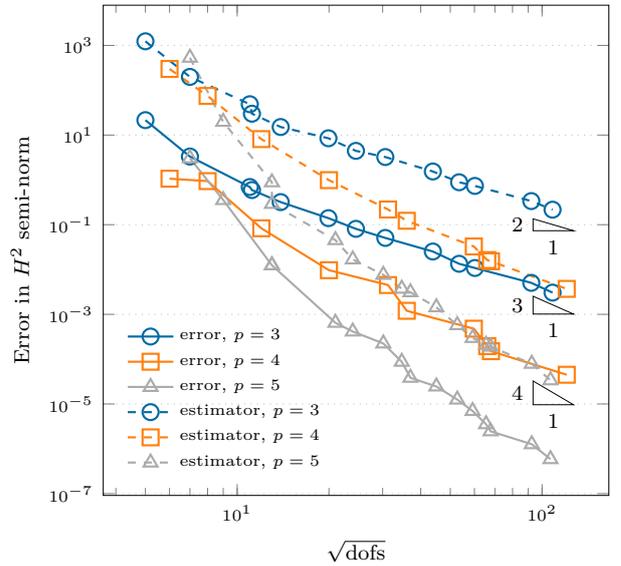}
		\caption{Residual-based estimator.}
	\end{subfigure}
	\caption{Study of the convergence of the error measured in the $H^2$ semi-norm, bubble and residual-based estimators against the square root of the number of dofs for different $p$ on a smooth plate problem ($u_{ex}=\sin(2 \pi x) \sin(2 \pi y)$). Adaptive refinement based on the maximum strategy ($\gamma = 0.5$).}
\label{fig:convergence_smooth_MS}
\end{figure}

\begin{figure}
	\centering
	\begin{subfigure}[t]{0.45\textwidth}
		\centering
	\input{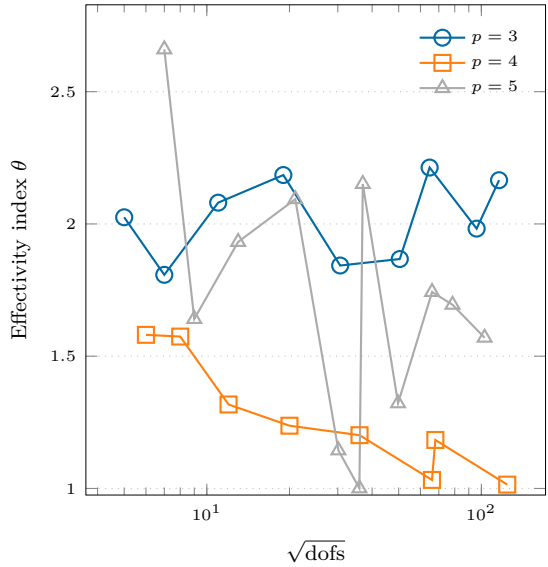}
		\caption{Bubble estimator.}
	\end{subfigure}
	\hfill
	\begin{subfigure}[t]{0.45\textwidth}
		\centering
	\input{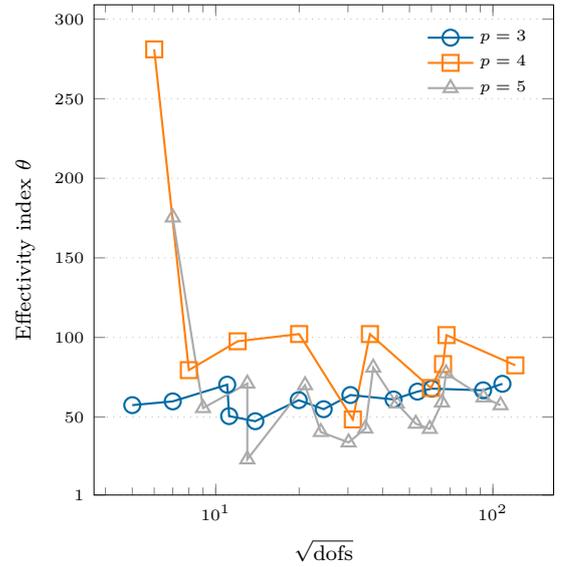}
		\caption{Residual-based estimator.}
	\end{subfigure}
	\caption{Effectivity index $\theta$ for the bubble and residual-based estimator against the square root of the dofs for different $p=3,4,5$ on the plate problem with smooth solution. Notice the difference of two orders of magnitude in the scale used for the y-axis.}
\label{fig:effectivity_index}
\end{figure}

Finally, we run the same example letting the bubble and residual-based estimators drive automatically the adaptive simulation. The results are reported in~\Cref{fig:convergence_smooth_MS}, where the error in energy norm and the estimator are plotted against the square root of the number of dofs. We observe analogous results to the uniform refinement case, where for $p=3,4,5$ the expected asymptotic rates of convergence are achieved. Then, let us define the effectivity index $\theta$ as:
\begin{align*}
\theta = \frac{ \sqrt{\sum_{\epsilon \in Q} \eta_{\epsilon}^2}}{\left\lVert u_{ex} - u_h \right\rVert_{E(\Omega)}} \, ,
\end{align*} 
where the optimal value of $\theta$ is $1$. Also in the adaptive example, the bubble estimator provides a better estimate of the real error compared to a classical residual-based estimator, yielding an effectivity index much closer to the optimal value, as shown in~\Cref{fig:effectivity_index}. This statement holds in an analogous way for all the presented examples and shows numerically the efficiency of the proposed method.

%

\FloatBarrier

\subsubsection{Singular solution on a square plate}

\begin{figure}
	\begin{subfigure}[t]{0.475\textwidth}
		\centering
	\includegraphics[scale=1.0]{pictures/numericalExamples/scordelis_geometry/plate.pdf}
	\caption{Geometry description of the problem.}
	\end{subfigure}
	\hfill
	\begin{subfigure}[t]{0.475\textwidth}
		\centering
	\includegraphics[scale=1.0]{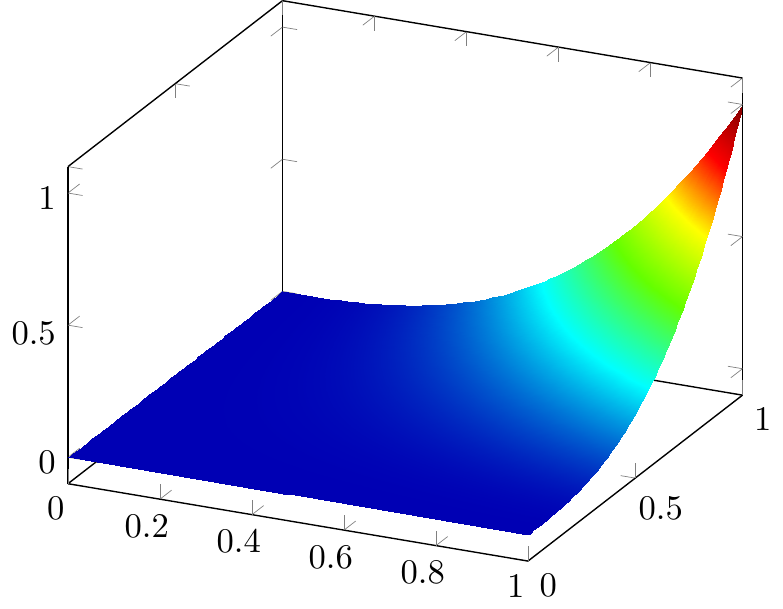}
	\caption{Exact solution $u_{ex} = x^\alpha y^\beta$ with $\alpha = \beta = 2.8$.}
	\end{subfigure}
	\caption{Geometry, physical parameters and exact solution contour of the singular square plate example.}
	\label{fig:plate_singular_geo}
\end{figure}

\begin{figure}[!h]
	\centering
		\input{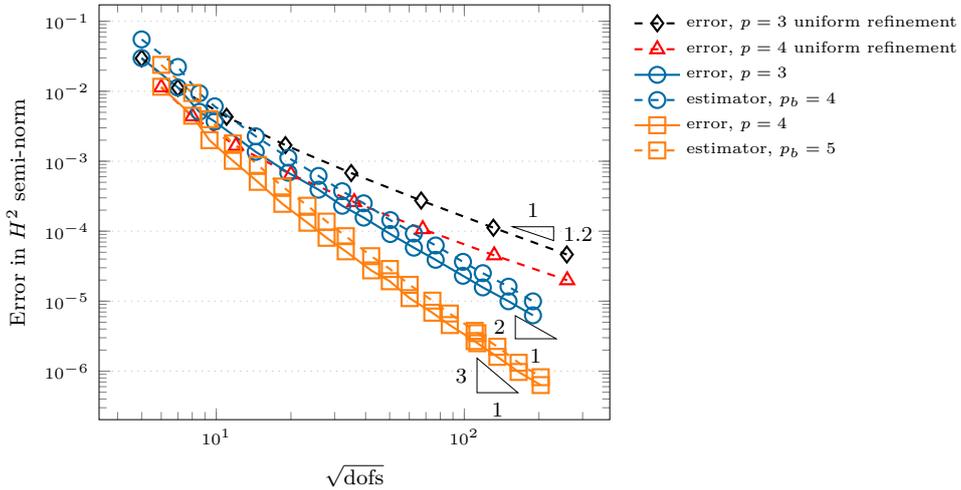}
	\caption{Study of the convergence of the error measured in the $H^2$ semi-norm and bubbles estimator against the square root of the number of dofs for different $p$ (and associated $p_b$ for the bubble space) on the square plate problem with singularity ($u_{ex}=x^\alpha y^\beta$ with $\alpha = \beta = 2.8$). Adaptive refinement based on the maximum strategy ($\gamma = 0.5$).}
	\label{fig:convergence_bubbles_singular_MS}
\end{figure}

\begin{figure}
		\centering
		\begin{subfigure}[t]{0.305\textwidth}
			\centering
			\includegraphics[width=1.2\textwidth]{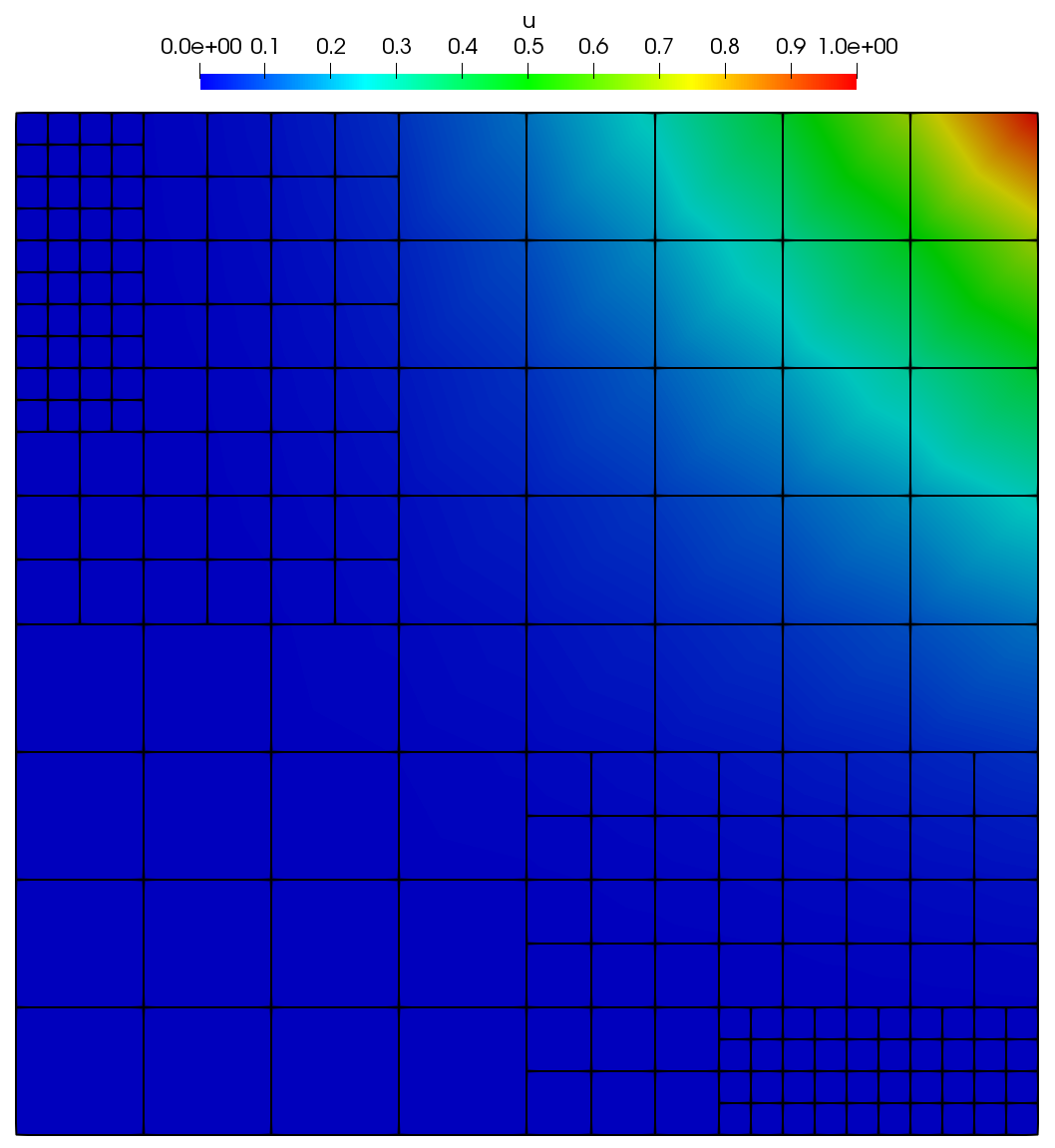}
			\caption{Solution $u$ and mesh at iteration $k = 3$.}
		\end{subfigure}
		\hspace{2.0cm}
		\begin{subfigure}[t]{0.305\textwidth}
			\centering
			\includegraphics[width=1.2\textwidth]{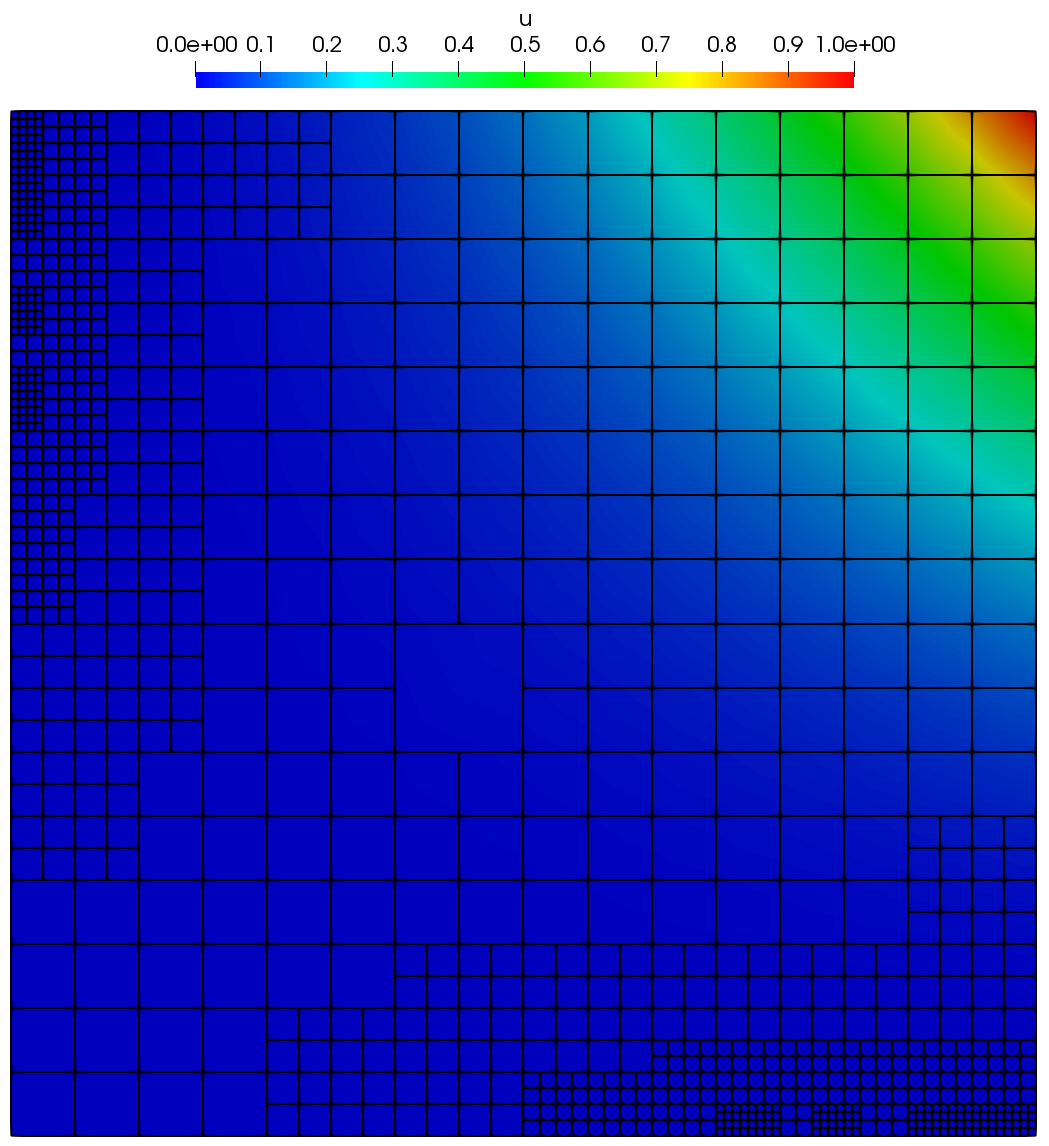}
			\caption{Solution $u$ and mesh at iteration $k = 5$.}
		\end{subfigure}
		\begin{subfigure}[t]{0.305\textwidth}
			\centering
			\includegraphics[width=1.2\textwidth]{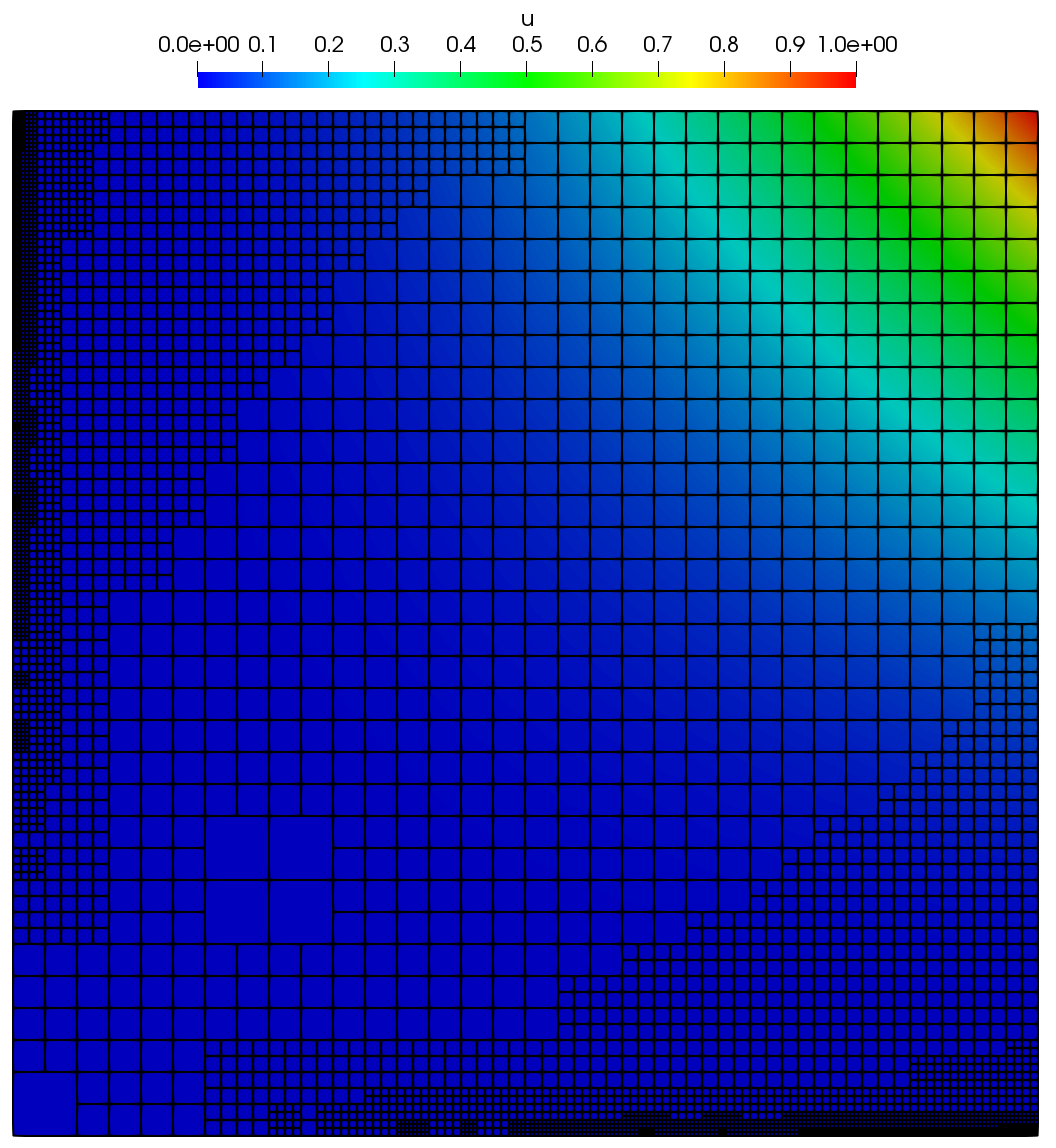}
			\caption{Solution $u$ and mesh at iteration $k = 8$.}
		\end{subfigure}
		\hspace{2.0cm}
		\begin{subfigure}[t]{0.305\textwidth}
			\centering
			\includegraphics[width=1.2\textwidth]{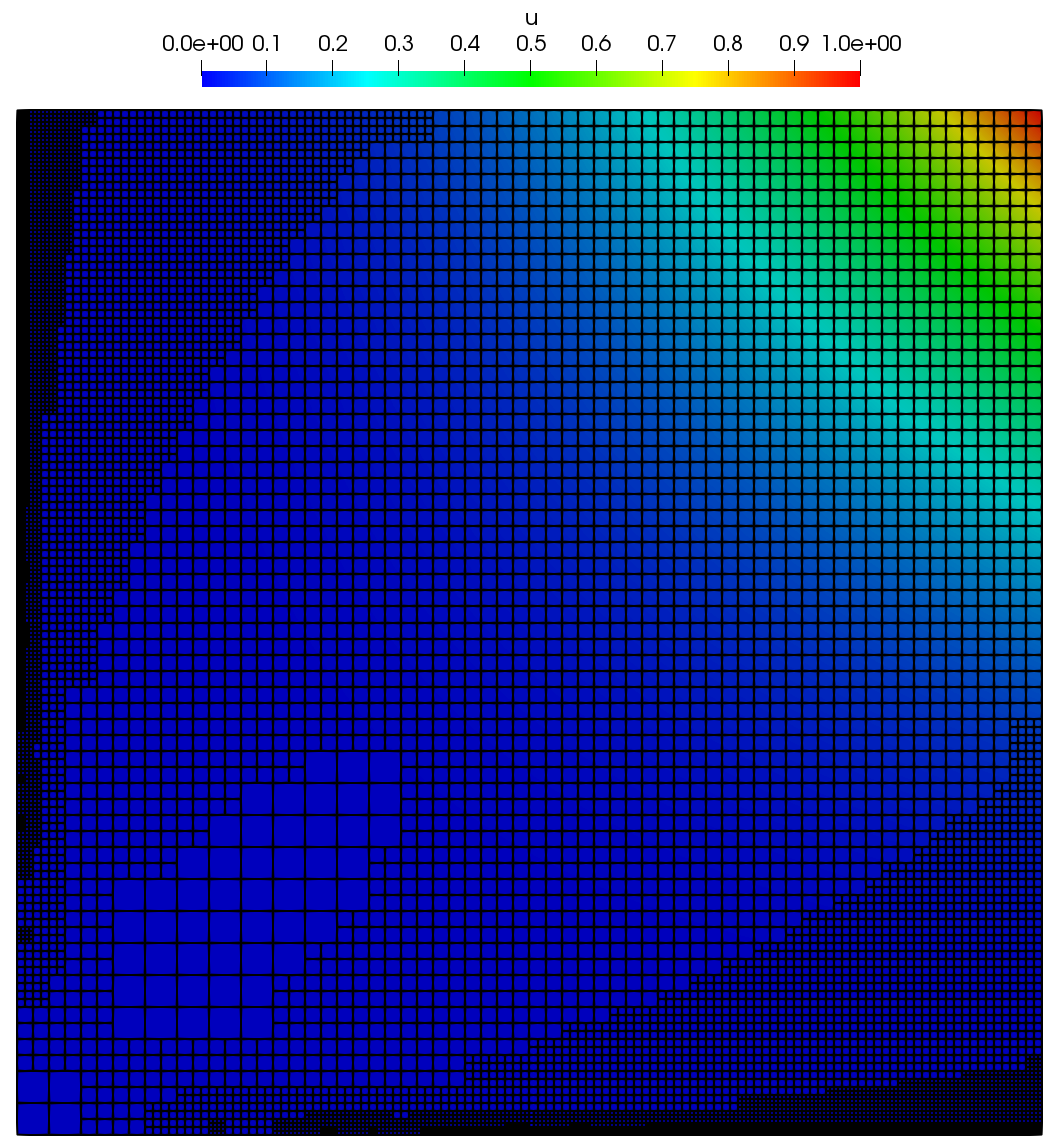}
			\caption{Solution $u$ and mesh at iteration $k = 10$.}
		\end{subfigure}
		\vfill
		\begin{subfigure}[t]{\textwidth}
			\centering
			\includegraphics[width=0.7\textwidth]{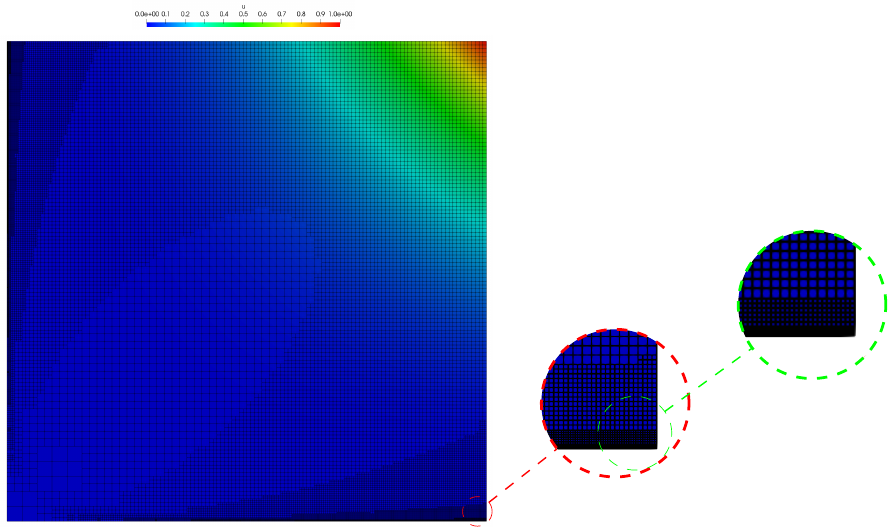}
			\caption{Solution $u$ and mesh at iteration $k = 14$ and zoom on the bottom right corner. Notice the very fine level of refinement achieved close to the singularity.}
		\end{subfigure}
			\caption{Mesh and solution contours at different steps of the adaptive loop driven by the bubble error estimator for the square plate with singular solution. Solution obtained using hierarchical B-splines of degree $p=3$.}
		\label{fig:mesh_singular_plate}
\end{figure}

In the next example, we consider again the computational domain to be a unit square $\Omega = [0,1]^2$, see Figure~\ref{fig:plate_singular_geo}. 
However, this time the manufactured solution $u_{ex} = x^\alpha y^\beta$ with $\alpha = \beta = 2.8$ is constructed such that a singularity is present along the bottom and the left edges of the plate. In particular, it holds $u_{ex} \in H^3(\Omega) \setminus H^4(\Omega)$. 
The applied load $g$ is again constructed such that it fulfills the manufactured solution $u_{ex}$ and it is given as follows:
\begin{align*}
g =& \, \alpha (\alpha - 1)(\alpha - 2)(\alpha - 3) x^{(\alpha - 4)} y^\beta + \beta(\beta - 1)(\beta - 2)(\beta - 3) x^\alpha y^{(\beta - 4)} + \\
& \,2 \left[ \alpha \beta(\alpha - 1)(\beta - 1) x^{(\alpha - 2)} y^{(\beta - 2)}\right]. 
\end{align*}
The boundary conditions are also constructed from the exact solution as:
\begin{alignat*}{2}
u &= u_{ex} \qquad &&\text{on} \quad \partial \Omega \quad , \\
\nu D \varDelta u + (1 - \nu) D \, \boldsymbol{d} \cdot (\nabla \nabla u)\boldsymbol{d} &= M_{ex} \qquad &&\text{on} \quad \partial \Omega 
\end{alignat*}
where $M_{ex}$ denotes the exact bending moment. The reduction in regularity of the solution limits the rate of convergence in case of uniform refinement, even for increasing $p$. From \citep{Ciarlet2002,Quarteroni2008}, the result on a priori convergence for uniform refinement is in this case a more general version of~\eqref{eq:apriori_smooth} and reads:
\begin{align*}
\left\lVert u_{ex} - u_h \right\rVert_{H^s(\Omega)} \leq C h^{l - s} \left\lVert u_{ex} \right\rVert_{H^{r}(\Omega)} \qquad r > s \, ,
\end{align*}
where $r$ represents the regularity of the exact solution (for this example we have $r \in \left] 3,4 \right[$) and $l$ is defined as $l = \min{(r,p+1)}.$
This effect can be clearly seen in~\Cref{fig:convergence_bubbles_singular_MS} for $p=3,4$. \Cref{fig:convergence_bubbles_singular_MS} also shows the results for an adaptive simulation driven by the bubble estimator, where optimal rates of convergence are recovered by our method and a significant increase in accuracy per degree-of-freedom is achieved. Additionally, the obtained mesh at different iterations of the adaptive algorithm are depicted in~\Cref{fig:mesh_singular_plate}, where it can be seen that the singularities are accurately detected and resolved by the estimator.

\subsubsection{Point load on a rectangular plate}

\begin{figure}[!h]
	\begin{subfigure}[t]{0.495\textwidth}
	\centering
		\input{\graphDir/convergence_bilaplacian_pointLoad_adm.tex}
	\caption{Study of the convergence of the displacement under a point load on the square plate example, adaptive refinement vs. uniform refinement using hierarchical B-splines of degree $p=3$.}
	\label{fig:pointLoad_displacement}
	\end{subfigure}
	\begin{subfigure}[t]{0.495\textwidth}
		\centering
		\includegraphics[width=0.9\textwidth]{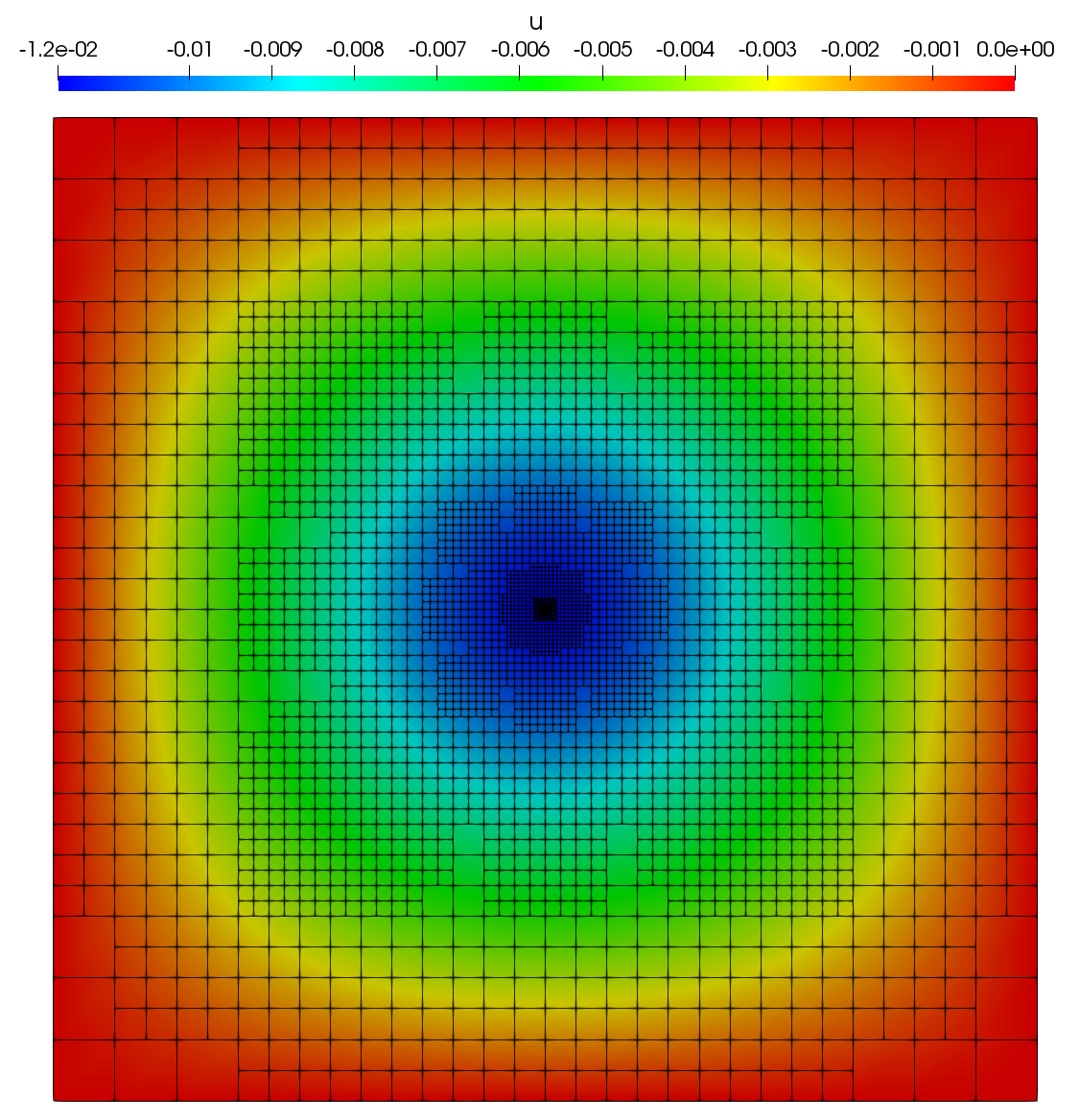}
		\caption{Mesh at iteration $k = 9$ and corresponding solution obtained by adaptive refinement driven by the bubble estimator using hierarchical B-splines of degree $p=3$.}
		\label{fig:pointLoad_mesh}
	\end{subfigure}
	\caption{Convergence study and numerical solution for the square plate problem subjected to a point load applied at the center of the plate.}	
\end{figure}

In this example, we test the performance of the proposed error estimator when a point load is applied as an external force. We perform again our computation on a unit square $\Omega = [0,1]^2$ and we suppose that the plate is simply supported on the entire boundary $\partial \Omega$. Additionally, we apply the external point load at the center of the structure. There exists an analytical solution for the deflection under the load, given as~\citep{reddy2006theory}: 
\begin{align}
u_{ex} = \frac{4 g L^2}{D \pi^4} \sum_{n=1,3,\ldots}^\infty \, \sum_{m=1,3,\ldots}^\infty \frac{1}{\left(m^2 + n^2\right)^2} \, ,
\label{eq:exact_point_load}
\end{align}
where $g$ represents the applied external force, $L$ is the length of the plate and $D$ denotes its flexural stiffness. We set $g=-1$, $L=1$ and the physical parameters such that $D=1$. Computing \eqref{eq:exact_point_load} with an adequate number of terms for the double Fourier series yields a reference value of $u_{ex} = -0.011600839735872 \ldots$ for the deflection. In~\Cref{fig:pointLoad_displacement} the convergence of the normalized deflection is plotted against the number of degrees-of-freedom, where the former is defined as $\vert 1 -u_h/u_{ex} \vert$. Moreover, the obtained solution and hierarchical mesh obtained after $k=9$ steps of the adaptive algorithm are depicted in~\Cref{fig:pointLoad_mesh}. Here, the advantages of using an adaptive scheme are clearly highlighted, in terms of efficiency and accuracy per degree-of-freedom. For instance, at around 2500 dofs the adaptive strategy is already two orders of magnitude more accurate compare to the results obtained by uniform refinement. 

\begin{remark}[Regularization of the Dirac delta]
It is worth noting that a point load is modeled as a \textit{Dirac delta}, which rigorously speaking is a distribution. Therefore, for the evaluation of the residual in a strong sense, needed in the residual-based estimator, we cannot directly use it inside our computations but instead we must regularize it. In our example, this is achieved with a steep Gaussian function. A considerable amount of literature has been written on how to regularize and correctly integrate the Dirac delta according to the corresponding application (e.g. we refer to \citep{Tornberg2004} and references therein) and although our approach might not be optimal we feel it suits the purpose of our test case while at the same time it proves to be easily implementable. 
\end{remark}

\FloatBarrier

\subsection{Kirchhoff-Love shell}

In this Section we demonstrate the applicability of the bubble error estimator to the Kirchhoff-Love formulation. Again, in all our numerical experiments the estimator performs well in driving adaptive computations. Once more, it is worth highlighting that with the proposed method we avoid the evaluation of higher order derivatives (this time defined on a manifold) which turns out to be particularly challenging and computationally involved in the framework of Kirchhoff-Love shells, see for instance \citep{Maurin2018}.

\subsubsection{Gravity load on a pinned roof}
\begin{figure}
	\centering
	\includegraphics[scale=0.9]{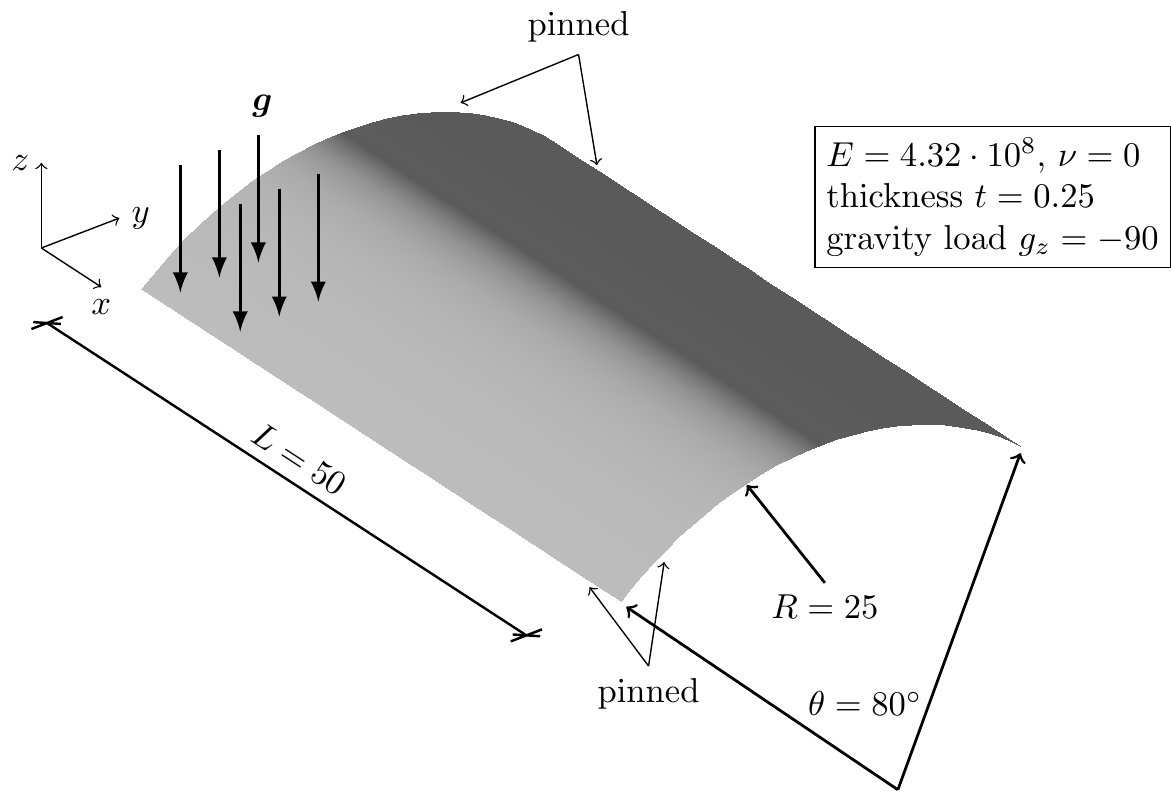}
	\caption{Geometry description and physical parameters of the pinned roof example (only part of the load is depicted for visualization purposes).}
	\label{fig:roof_pinned_geo}
\end{figure}

In this example, we study the behavior of a structure subjected to its self-weight where pinned supports are applied to the entire boundary. The geometry, boundary conditions and physical parameters are given in~\Cref{fig:roof_pinned_geo}, where the problem setup (except for the boundary conditions) is taken from the classical Scordelis-Lo benchmark~\citep{Timoshenko1959}. The boundary conditions are modified in order to avoid numerical issues during the computation of the reference solution and corresponding error, since the original problem definition of the Scordelis-Lo roof is not well-posed, see for instance~\citep{Kiendl2009}.  
To the best of the authors' knowledge, there is no closed-form global solution available in the literature for the problem at hand. Therefore, all the convergence studies presented in the following assess the behavior against a reference solution $\boldsymbol{u}_h^{\text{ref}}$, which was computed using B-splines of degree $p=8$ on a fine uniform mesh with $\simeq 200 \, 000$ degrees-of-freedom. In particular, we define an approximation of the error in the energy norm as:
\begin{align*}
\vert \vert \tilde{e} \vert \vert_{E(\Omega)} = \sqrt{a(\boldsymbol{u}_h - \boldsymbol{u}_h^{\text{ref}},\boldsymbol{u}_h - \boldsymbol{u}_h^{\text{ref}})} \quad ,
\end{align*} 
where the numerical integration of this quantity is performed on the fine mesh. 

\begin{figure}
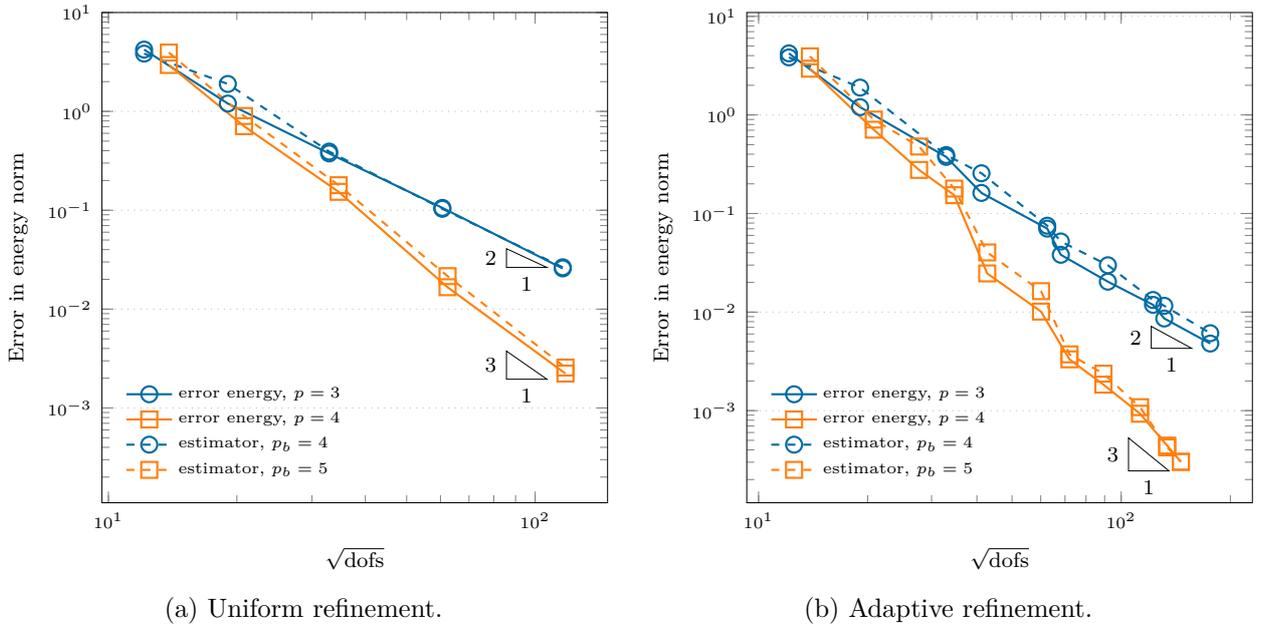

	\centering
	\begin{subfigure}[t]{0.495\textwidth}
		\centering
	\input{\graphDir/convergence_scordelisLo_smooth_bubbles_GR.tex}
		\caption{Uniform refinement.}
	\end{subfigure}
	\hfill
	\begin{subfigure}[t]{0.495\textwidth}
		\centering
	\input{\graphDir/convergence_scordelisLo_smooth_bubbles_MS.tex}
		\caption{Adaptive refinement.}
	\end{subfigure}
	\caption{Study of the convergence of the error in the energy norm for the bubble estimator employing hierarchical B-splines of degree $p=3,4$ on the pinned roof example. Adaptive refinement based on the maximum strategy ($\gamma = 0.5$).}
\label{fig:convergence_roof_pinned}
\end{figure}

In~\Cref{fig:convergence_roof_pinned} the convergence behavior of $\vert \vert \tilde{e} \vert \vert$ is depicted for the case of uniform refinement and adaptive refinement driven by the bubble estimator, against the number of dofs. We notice that since the true solution of the problem is regular enough, we obtained the optimal asymptotic rate of convergence, both for the error in energy norm and the estimator.
This example is thought as a first assessment of the performance of the proposed estimator on shell geometries. We remark that in standard shell benchmarks adopted in the literature, convergence is usually tested against a reference value for the displacement in some points of interest of the structure. While this information is relevant in many engineering applications, we feel that a deeper look into the behavior of derived global quantities, like for instance the error in energy norm used here, can be useful and mathematically more rigorous.

\subsubsection{Point load on the Scordelis-Lo roof}

\begin{figure}[!h]
	\centering
	\begin{subfigure}[t]{0.475\textwidth}
	\includegraphics[scale=0.75]{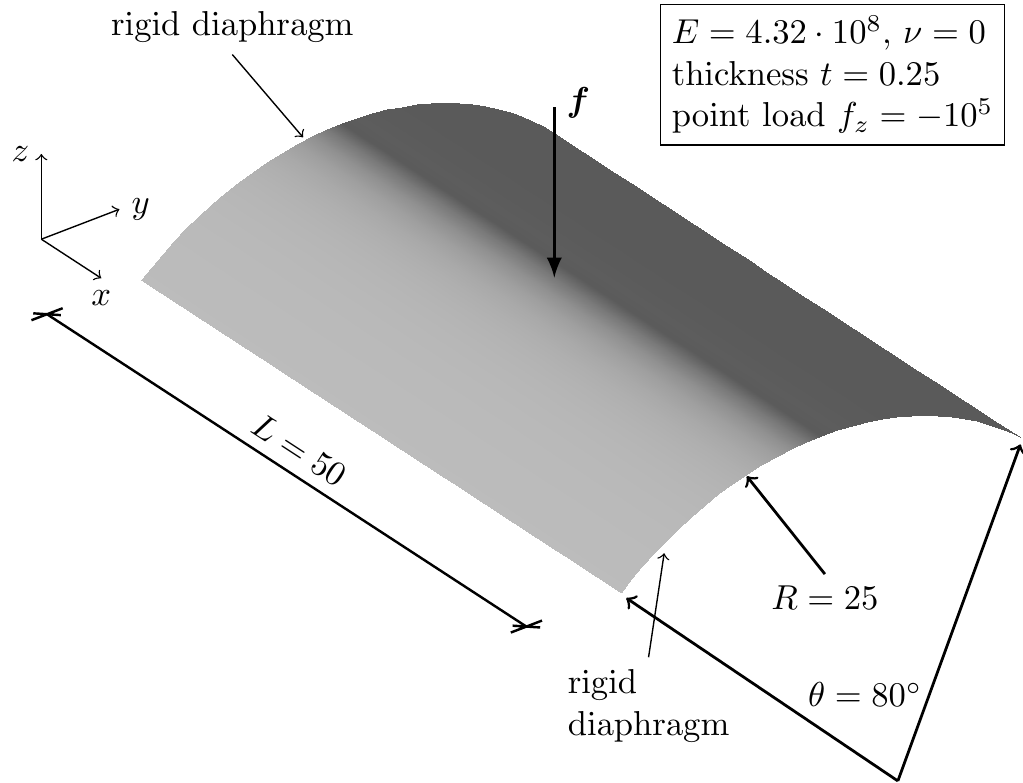}
	\caption{Geometry description and physical parameters of the Scordelis-Lo roof subjected to a point load.}
	\label{fig:scordelis_pointLoad_geo}	
	\end{subfigure}
	\begin{subfigure}[t]{0.475\textwidth}
		\input{\graphDir/convergence_scordelisLo_pointLoad_log.tex}
	\caption{Study of the convergence of the displacement on the shell under a point load, uniform refinement against adaptive refinement (maximum strategy, $\gamma = 0.5$) using hierarchical \mbox{B-splines} of degree $p=3$.}
	\label{fig:convergence_scordelisLo_bubbles_GR_meshSize}
	\end{subfigure}
	\caption{Problem setup and convergence plot for the Scordelis-Lo roof example under a point load.}
\end{figure}

\begin{figure}
		\centering
		\begin{subfigure}[t]{0.325\textwidth}
			\centering
			\includegraphics[width=\textwidth]{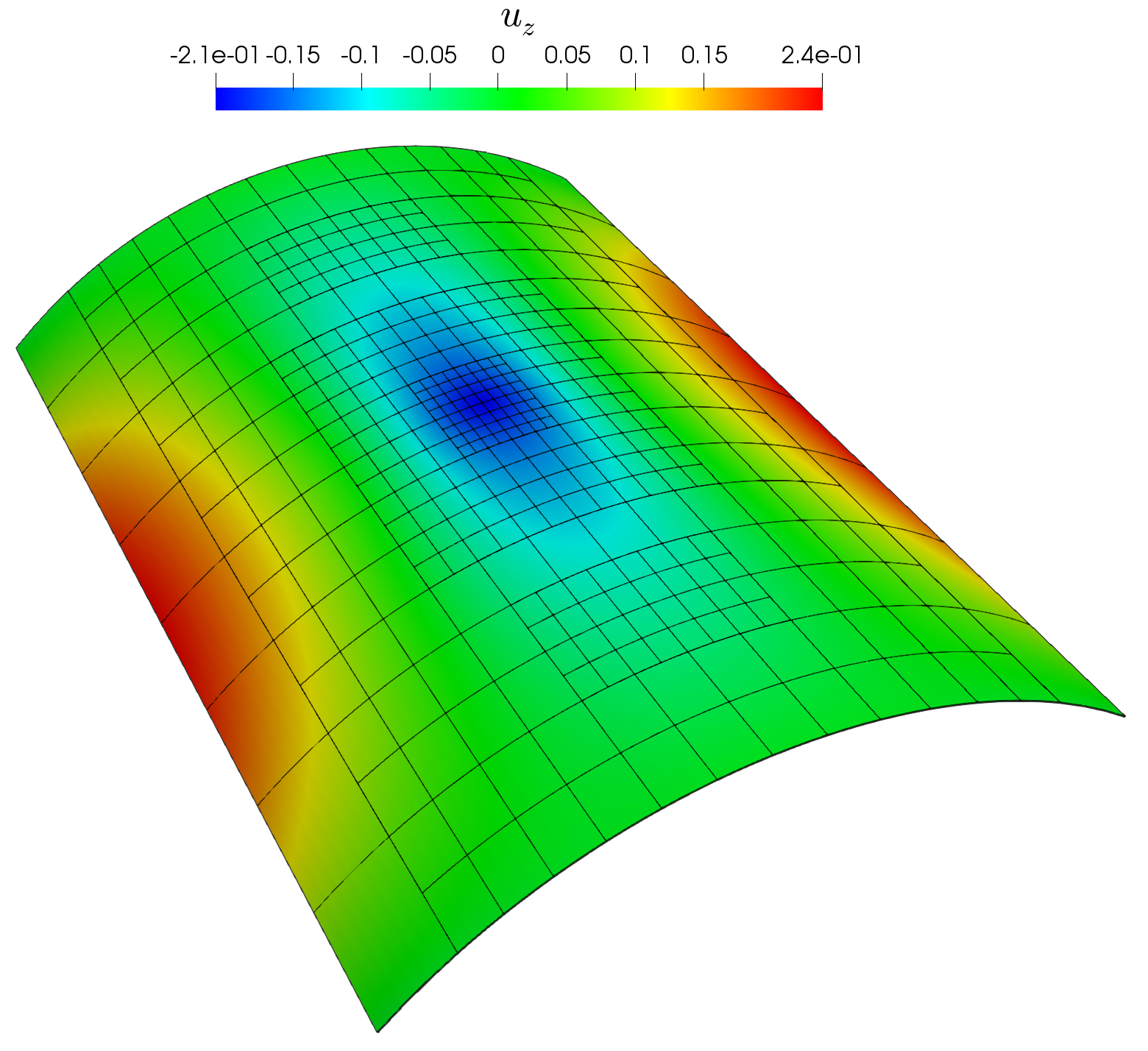}
			\caption{Solution $u_z$ and mesh at iteration $k = 3$.}
		\end{subfigure}
		\begin{subfigure}[t]{0.325\textwidth}
			\centering
			\includegraphics[width=\textwidth]{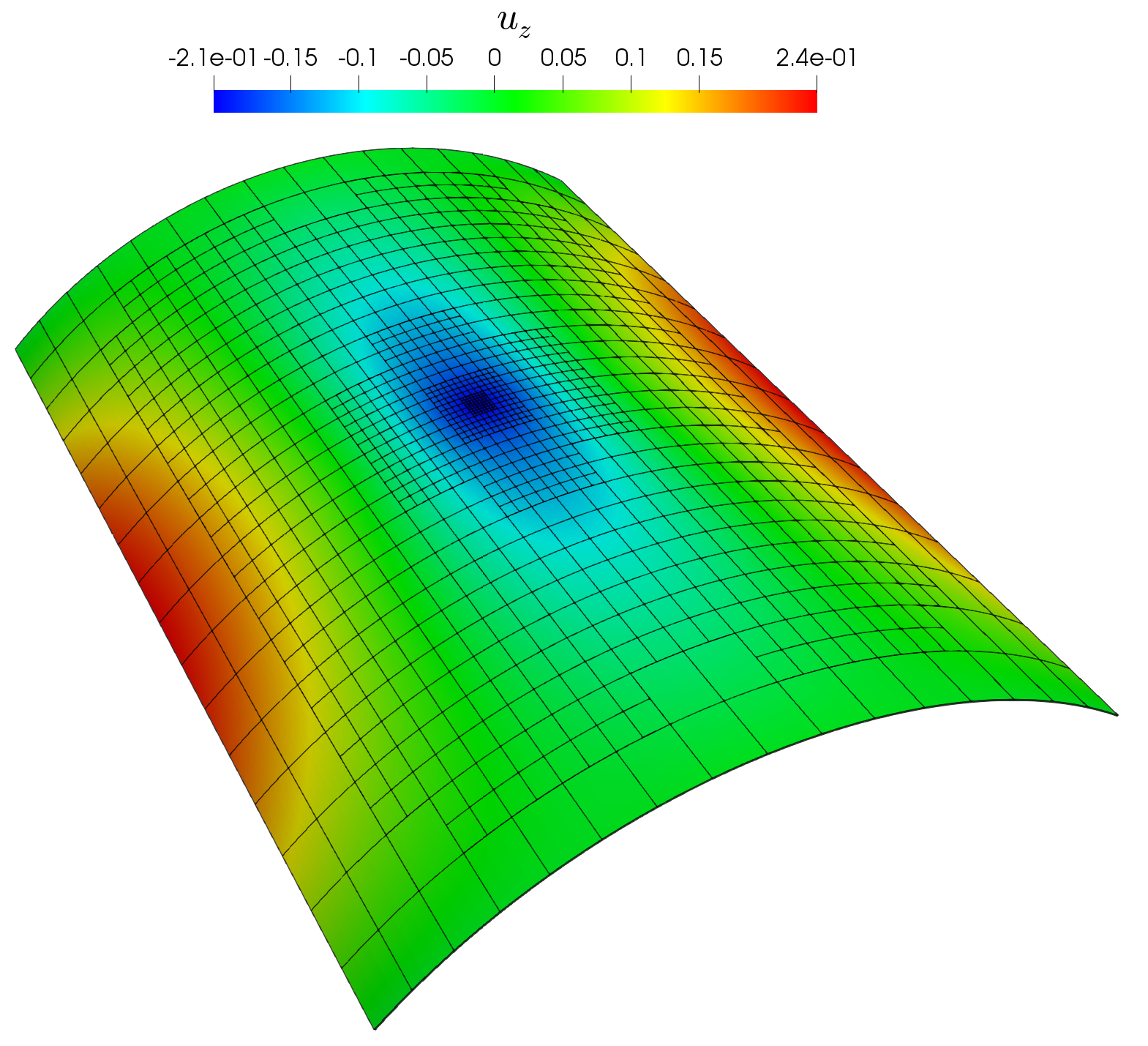}
			\caption{Solution $u_z$ and mesh at iteration $k = 5$.}
		\end{subfigure}
		\begin{subfigure}[t]{0.325\textwidth}
			\centering
			\includegraphics[width=\textwidth]{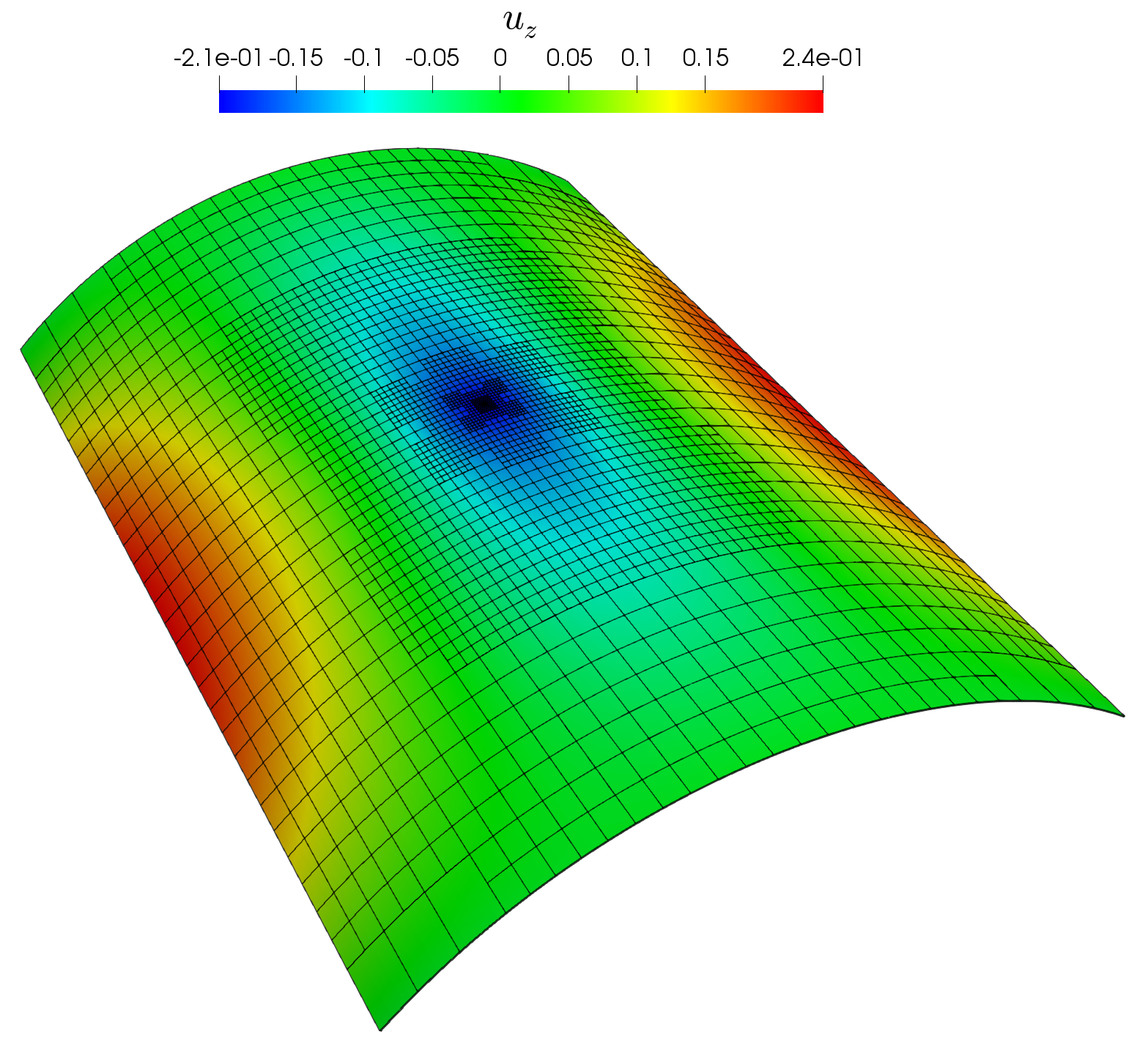}
			\caption{Solution $u_z$ and mesh at iteration $k = 9$.}
		\end{subfigure}
		\begin{subfigure}[t]{0.325\textwidth}
			\centering
			\includegraphics[width=\textwidth]{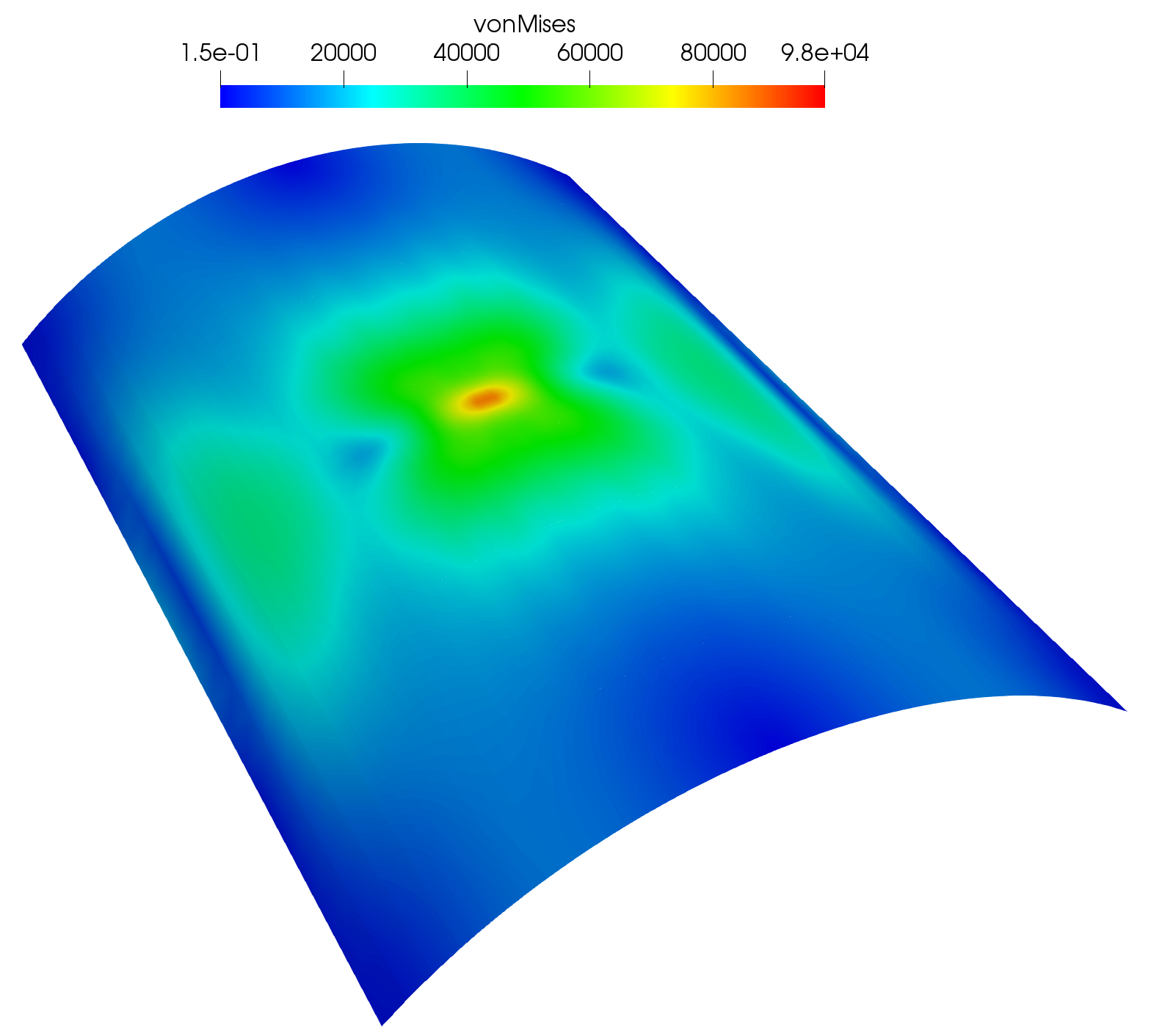}
			\caption{Von Mises at iteration $k = 3$.}
		\end{subfigure}
		\begin{subfigure}[t]{0.325\textwidth}
			\centering
			\includegraphics[width=\textwidth]{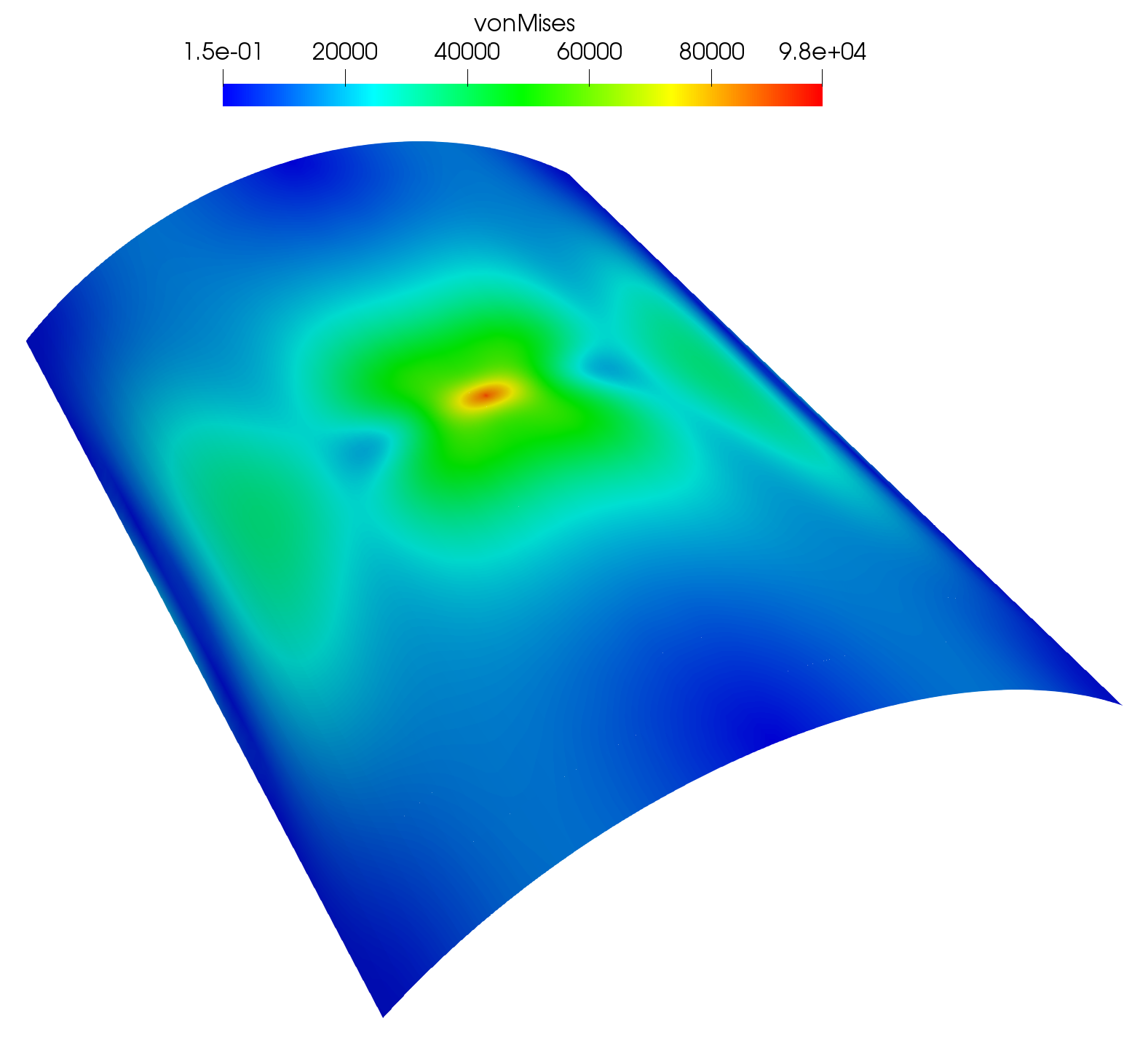}
			\caption{Von Mises at iteration $k = 5$.}
		\end{subfigure}
		\begin{subfigure}[t]{0.325\textwidth}
			\centering
			\includegraphics[width=\textwidth]{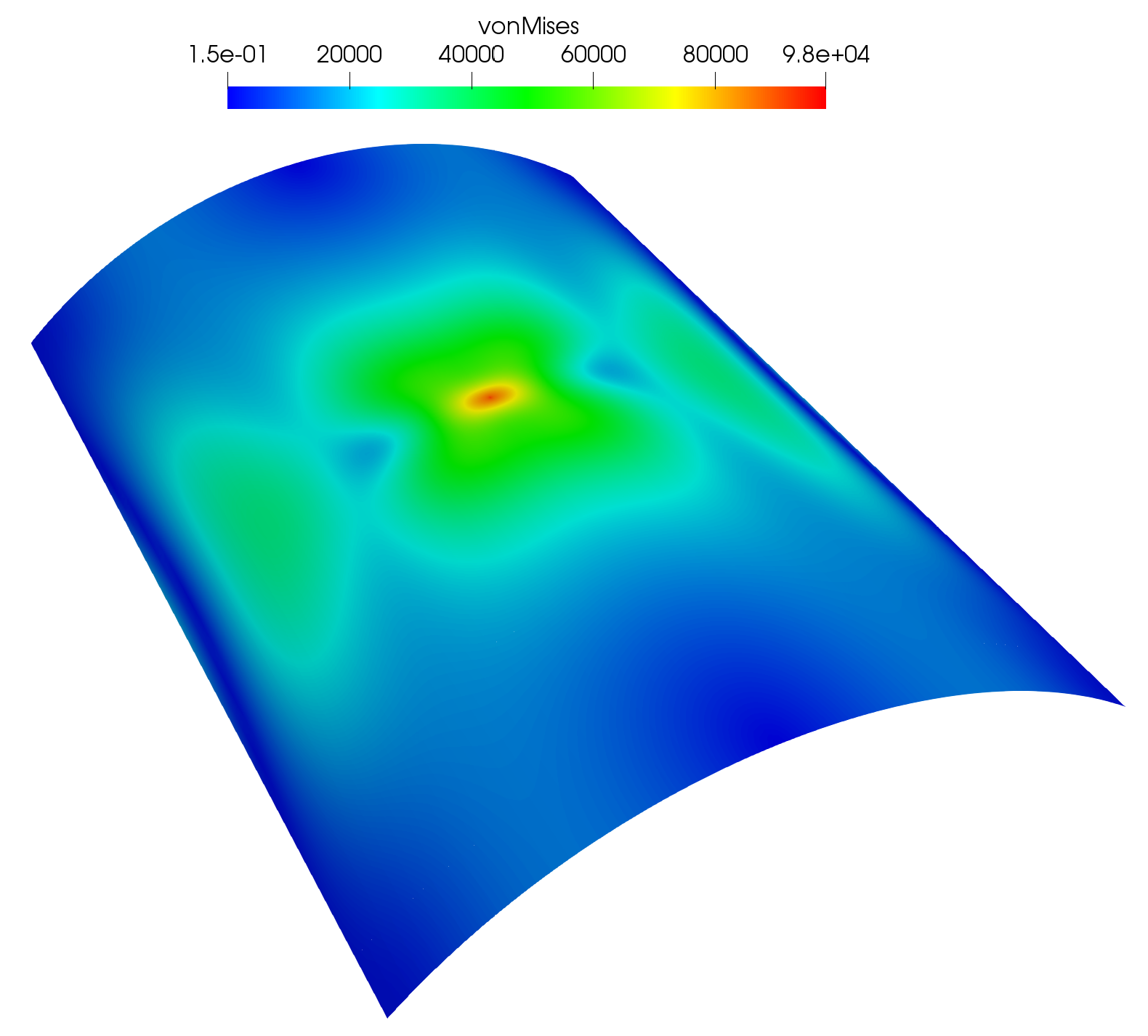}
			\caption{Von Mises at iteration $k = 9$.}
		\end{subfigure}		
		\begin{subfigure}[t]{0.75\textwidth}
			\centering
			\includegraphics[width=1.1\textwidth]{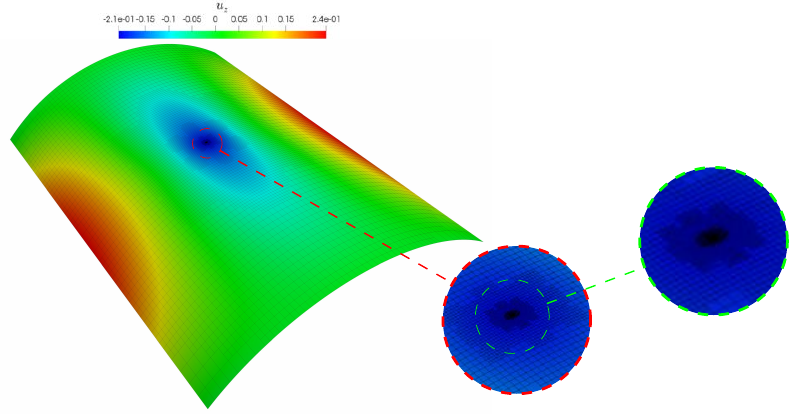}
			\caption{Solution $u_z$ and mesh at iteration $k = 11$.}
		\end{subfigure}
			\caption{Mesh, solution $u_z$ and Von Mises stress at different steps of the adaptive loop driven by the bubble error estimator for the Scordelis-Lo roof subjected to a point load, solution obtained employing hierarchical B-splines of degree $p=3$.}
		\label{fig:scordelis_mesh_and_vonMises}
\end{figure}

The next example is meant to demonstrate once again the higher accuracy per-degree-of-freedom achievable using local refinement. The geometrical setup and physical parameters are taken as defined in the Scordelis-Lo benchmark and are given in~\Cref{fig:scordelis_pointLoad_geo}, where we change the boundary conditions. Indeed, as external loading, we apply a point load of magnitude $10^5$ in the middle of the structure, directed in the negative vertical direction.
The convergence behavior of the displacement under the point load is depicted in~\Cref{fig:convergence_scordelisLo_bubbles_GR_meshSize}, where a reference value of $-0.206794699852788 \ldots$ has been obtained with an overkill solution, which was computed using B-splines of degree $p=4$ on a fine uniform mesh with $\simeq 200 \, 000$ degrees-of-freedom.
It can be seen that the solution obtained with an adaptive simulation based on the bubble estimator is several orders of magnitude more accurate compared to performing uniform refinement of the patch, for the same number of dofs.
Finally, in~\Cref{fig:scordelis_mesh_and_vonMises}, the obtained displacement in the z-direction, hierarchical mesh and Von Mises stress distribution are presented at different steps of the adaptive loop where we remark that once again the estimator properly captures and resolves the sharp features of the solution.

\section{Conclusions} \label{sec:conclusions}
We have introduced a novel \textit{a posteriori} error estimator in the context of linear fourth-order partial differential equations, namely we have successfully applied it to adaptive simulation, based on the use of hierarchical B-splines, of isogeometric Kirchhoff plates and Kirchhoff-Love shells. This can be seen as an extension of similar approaches present in the literature, e.g. we refer to \citep{Bank1993,Vuong2011}.
The evaluation of the estimator is based on the solution of an additional, residual-like variational problem formulated in the so-called bubble space, which is composed of Bernstein polynomials defined locally on active elements. We remark that, thanks to this choice of the aforementioned space, the resulting linear system is in general small, block-diagonal and easily-invertible.
Moreover, this method is suitable for parallelization and straightforward to implement into existing isogeometric codes. More importantly, it is computationally cheap compared to classical residual-type estimators since it avoids the computation of the residual in a strong sense. On one hand, this is a major advantage particularly for Kirchhoff-Love shells, since the evaluation of covariant derivatives is a tedious task, becomes very quickly computationally expensive and, for all practical purposes, almost intractable from a numerical standpoint. On the other hand, with this technique we also avoid the computation of integral terms involving the evaluation of the jump of the derivatives across element boundaries. Finally, we have also observed that in all our experiments the proposed estimator acts as an excellent approximation of the true error. 

To conclude, we have numerically demonstrated the applicability and robustness of the proposed error estimator to steer an adaptive simulation of a wide range of problems of engineering relevance, where superior efficiency and accuracy per-degree-of-freedom have been achieved thanks to the local refinement capabilities of hierarchical B-splines.

\section*{Acknowledgements} 
The authors would like to thank Prof.Alessandro Reali and Dr.Rafael V\'{a}zquez for the fruitful discussions on the subject of this paper and Dr.Rafael V\'{a}zquez for his help with GeoPDEs. The authors also gratefully acknowledge the support of the European Research Council, via the ERC AdG project CHANGE n.694515.

\bibliographystyle{plainnatnourl}
\bibliography{library.bib}

\end{document}